\documentclass[]{ajour}
\usepackage{amsmath, amssymb, amsthm, conformal, mystyle}

\renewcommand{\theenumi}{{\rm(\alph{enumi})}}

\renewcommand\^[1]{\widehat{#1}}

\def\c{{\mathsf c}}

\begin{document}

%\bibliographystyle{plain}

%%%%% To be entered at Academic Press: =====>>
%
% Uncomment line below only when doing final typesetting,
%\finaltypesetting
% \journame{}
% \articlenumber{}
% \yearofpublication{}
% \volume{}
% \cccline{}
% \received{}
% \revised{}
% \accepted{}

% communication line, use: \commline{Communicated by...}
% \commline{Communicated by... }

\authorrunninghead{Michael Roitman}
\titlerunninghead{On twisted representations of vertex algebras}

%% To set particular page number:
%\setcounter{page}{261} %%

%% <<== End of commands to be entered at Academic Press

\title{On twisted representations of vertex algebras\thanks{This
research was conducted by the author for the Clay Mathematics Institute}}

\author{Michael Roitman}

\affil{Department of Mathematics \\
University of Michigan \\
Ann Arbor, Michigan  48109-1109
}
\email{roitman@umich.edu}

\abstract{
In this paper we develop a formalism for working with representations
of vertex and conformal algebras by generalized fields --- formal
power series involving non-integer powers of the variable.  
The main application of our technique is the construction of
a large family of representations for the vertex superalgebra
$\goth V_\Lambda$ corresponding to an integer lattice $\Lambda$. For
an automorphism $\^\sigma:\goth V_\Lambda\to\goth V_\Lambda$ coming
from a finite order automorphism $\sigma:\Lambda \to \Lambda$ we
find the conditions for existence of twisted modules of $\goth V_\Lambda$.
We show that the category of twisted representations of
$\goth V_\Lambda$  is semisimple with
finitely many isomorphism classes of simple objects.

}

%\date\today

%\subjclass{}

\begin{article}

\section*{Introduction}\label{sec:intro}
  
One of the most fascinating features of vertex algebras is their 
``sensitivity'' to finite order automorphisms. 
This is best illustrated by the construction of Moonshine representation of the
Monster simple group \cite{bor, flm}. 

In general, vertex
algebras are represented by formal power series, also called vertex operators, 
of the form  
%\vskip-10pt
$$
\alpha(z) = \sum_{n\in \Z} \alpha(n)\, z^{-n-1},\qquad \alpha(n) \in
\op{gl}(V),
$$
%\vskip-1pt\noindent
with coefficients in the algebra of linear operators on some space
$V$. In this case $V$ is a module over the vertex algebra. 
Having a finite order automorphism on a vertex algebra $\goth A$ means that
there is a grading $\goth A=\bigoplus_{\lambda\in \Gamma} \goth A^\lambda$
on $\goth A$ by a finite cyclic group $\Gamma$. Assume that $\Gamma\subset
\C/\Z$, then it is tempting to represent the vertex algebra $\goth A$ by
the so-called {\it twisted} or {\it generalized} vertex operators, involving
non-integer powers of $z$, so 
that an element $\alpha \in \goth A^\lambda$ is represented by series of
the form 
\vskip-5pt
$$
\alpha(z) = \sum_{n\equiv \lambda \op{mod} \Z} \alpha(n)\,
z^{-n-1},\qquad \alpha(n) \in \op{gl}(V). 
$$
This indeed could be done, and then $V$
becomes a {\it twisted module} over $\goth A$. In fact this construction
works for any group
$\Gamma\subset \C/\Z$. 

The idea of twisted realizations of vertex algebras goes back to the
work of Lepowsky and Wilson \cite{lepwil}, who introduced the
so-called twisted vertex operators. These operators
were systematically studied in \cite{flm,kp112,lep}. Twisted modules of
vertex algebras were
defined in \cite{dongtwisted, ffr}, see also 
\cite{dl2,kp112,litwisted}. 
Closely related is the
theory of orbifolds --- the invariant subalgebras of vertex algebras
under an action of a finite group of automorphisms. Twisted
representations of the vertex algebra yield ordinary (non-twisted)
representations of its orbifold, see e.g. \cite{bkt,dvvv}.

In this paper we develop a formalism for working with generalized
vertex operators
and use it to construct generalized representations of conformal and vertex
algebras in a very general setting. Then we consider two applications
of this techniques. 

First we study realizations of conformal algebras
by twisted formal series. Conformal algebras were  introduced by Kac
in the book \cite{kac2}, see also \cite{kac_fd,primc,freecv}.
They proved to be a valuable tool in studying
vertex algebras, the relation between the former and the latter is
somewhat like the relation between Lie and associative algebras. 

The second application is the construction of generalized
representations  of the vertex (super)algebra $\goth V_\Lambda$ 
corresponding to an integer lattice $\Lambda$. Our result is
similar to that of  Dong and Lepowsky
\cite{dl2}, but we use different
techniques and get a slightly more general construction. Lattice vertex algebras were
extensively studied in e.g. \cite{dong,dl,flm,kac2,lixu}. 
They play a very important role in
different areas of mathematics and physics, in particular the
Moonshine vertex algebra $V^\natural$, mentioned above, 
is closely related to the lattice vertex algebra of certain even
unimodular lattice of rank 24, called the Leech lattice. 
In fact one needs to consider a twisted representation of the vertex
algebra of the Leech lattice in order to construct the Moonshine
vertex algebra. 

An automorphism $\sigma:\Lambda\to \Lambda$ of the lattice can be
extended (in a non-unique way) to an automorphism $\^\sigma:\goth
V_\Lambda\to \goth V_\Lambda$ of the lattice vertex algebra. Our
construction of generalized $\goth V_\Lambda$-modules yields all
reasonable twisted modules of $\goth V_\Lambda$ corresponding
to the grading  by the root spaces of
$\^\sigma$. It turns out that sometimes there is no twisted
representation for either continuation  of $\sigma$. In this
case the representations of the orbifold vertex algebra 
$\goth V_\Lambda^{\^\sigma}$ do not come from the twisted
representation of $\goth V_\Lambda$. 

For an automorphism $\^\sigma:\goth V_\Lambda\to \goth V_\Lambda$
as above we define a category $\cal O_{\^\sigma}$ of twisted
representation, analogous to the category $\cal O$ of representations
of Kac-Moody Lie algebras. All reasonable twisted 
$\goth V_\Lambda$-modules, including those that satisfy a traditional
definition of a module over a vertex operator algebra as in e.g. 
\cite{ffr,flm}, belong to  $\cal O_{\^\sigma}$. We prove that
the category  $\cal O_{\^\sigma}$ is semisimple with finitely many
isomorphism classes of simple objects, that is, every module $V\in \cal O_{\^\sigma}$ is
decomposed into a direct sum of irreducible submodules, and there are
only finitely many irreducible modules, up to an isomorphism.
This result has been also 
obtained recently by Bakalov, Kac and Todorov \cite{bkt}. Some special
cases were known before, for example
the case when $\sigma=-1$ was studied by Dong and Nagatomo in \cite{dongtwisted,donag}, 
the case when $\Lambda$ is a simply-laced root lattice and $\sigma$ is
an element of the Weyl group of corresponding affine Kac-Moody algebra
was studied by Kac and Peterson \cite{kp112}.

\vfill%\smallskip
\noindent{\it Organization of the manuscript.}\quad
We start with giving formal 
definitions of conformal and vertex algebras. For more details the
reader can consult the books \cite{dl,flm,kac2}. Then, in 
\sec{locality}--\sec{modules} we discuss some properties of these algebras
in the context of generalized formal series. In \sec{products} we
derive a nice formula for the products of generalized series, which is
probably new. In \sec{series} we prove that conformal and vertex
algebras are exactly the algebraic structures formed by generalized
series with coefficients in a Lie algebra and generalized vertex operators
respectively. For vertex algebras this was proved by Li
\cite{litwisted} in a slightly less general setting and using different
methods. For the non-twisted case this result is well known
\cite{kac2,li}.

In \sec{2} we show how the approach developed in \sec{1} works for
conformal algebras. As in the non-twisted case, there is a universal
realization of a conformal algebra $\goth L$ with coefficients in a
certain Lie algebra $\cff_\Gamma \goth L$. This Lie algebra can be
constructed explicitly from $\goth L$. In \sec{affine} we illustrate
this by the example of an affine conformal algebra.  Similar ideas
appeared also in \cite{kac_fd}.

In \sec{3} we study generalized representations of lattice vertex
superalgebra $\goth V_\Lambda$. After some preliminary
information on representations of Heisenberg algebras (\sec{heisen})
and Fock spaces (\sec{fock}) we define in \sec{tvo} the twisted
vertex operators, first introduced by Lepowsky \cite{lep}. In 
\sec{lattice} we show that these operators generate a representation of the
lattice vertex algebra.

In \sec{lepwil} we show that we have in fact constructed all
reasonable generalized representations of lattice vertex algebras. The
argument uses an idea of Lepowsky and Wilson \cite{lepwil81}, which was
also used in \cite{dong,dongtwisted,flm,lixu}. 

Finally, in \sec{auto}--\sec{ss} we study the twisted representations
of $\goth V_\Lambda$, i.e. generalized representations, which are homogeneous with respect to the grading
induced by an automorphism $\^\sigma$ of $\goth V_\Lambda$. 
In \sec{auto} we find the conditions on generalized $\goth
V_\Lambda$-module, constructed in \sec{lattice}, to be $\^\sigma$-twisted,
while in \sec{ss} we define the category
$\cal O_{\^\sigma}$ of twisted $\goth V_\Lambda$-modules and 
prove that this category is semisimple with finitely many
isomorphism classes of simple objects. 

The key idea
is to show that the category $\cal O_{\^\sigma}$  is
equivalent to the category of graded representations of certain graded
associative algebra $A$, which turns out to be graded semisimple. This
idea is similar to the idea of Zhu algebra, introduced by Zhu \cite{zhu}
and then generalized for the twisted case by Dong, Li and Mason
\cite{dlm}.  However, our
algebra $A$ is quite different from the Zhu algebra of $\goth V_\Lambda$.

\section{Vertex algebras}\label{sec:1}
In this section we give an abstract definition of the main objects of
this paper --- conformal and
vertex algebras, and then show how they can be represented by generalized
formal power series. Conformal algebras were defined by Kac
\cite{kac2,kac_fd}, see also \cite{primc,freecv}. The first axiomatic
definition of vertex algebras is due to Borcherds \cite{bor}, see also 
\cite{fhl, flm, kac2}.
We recall here some basic properties of these
algebras and develop some techniques that will be used in \sec{2} and \sec{3}. 

All algebras and spaces are over a ground field $\Bbbk$ of
characteristic 0.

\subsection{Definitions of conformal and vertex algebras}
\label{sec:def}

A {\it Conformal algebra} 
is  a vector space $\goth L$ equipped with a sequence of
bilinear products $\ensquare n:\goth L\otimes\goth L \to \goth L, \ 
n\in \Z_+$, and a linear operator $D:\goth L\to \goth L$, such that the
following axioms hold for all $a,b,c\in\goth L$ and $n\in\Z_+$:
\begin{itemize}
\item[C1.](Locality)
$a\ensquare n b = 0$ for $n\gg0$.
\item[C2.]
$D(a\ensquare n b) = (Da) \ensquare n b + a \ensquare n (Db)= -n
a\ensquare{n-1} b+ a \ensquare n (Db)$.
\item[C3.](Quasisymmetry) 
\begin{equation*}
%$\displaystyle{
a\ensquare n b = - \sum_{i\ge0} (-1)^{n+i}\frac 1{i!} D^i
(b\ensquare{n+i}a).
%}$
\end{equation*}
\item[C4.](Conformal Jacoby identity)
\begin{equation*}%\label{fl:confjac}
(a \ensquare{n} b)\ensquare m c = 
\sum_{i=0}^n (-1)^i\binom{n}{i} 
\biggl(a\ensquare{n-i} (b \ensquare{m+i} c) -
b\ensquare{m+i} (a\ensquare{n-i} c)\biggr).
\end{equation*}
\end{itemize}

\vskip-10pt
\begin{Ex}
The Virasoro conformal algebra $\goth{Vir}$
is generated over $\Bbbk[D]$ by elements $\upsilon$ and $\c$, such that
$D\c=0$ (and therefore all products with $\c$ are 0 due to C2) and the
products of $\upsilon$ with itself are 
\begin{equation}\label{fl:vir}
\upsilon\ensquare 0 \upsilon = D\upsilon,\quad
\upsilon \ensquare 1 \upsilon = 2\upsilon,\quad
\upsilon\ensquare 3 \upsilon = \c.
\end{equation}
The rest of the products are 0.
\end{Ex}

A {\it vertex algebra} can be defined axiomatically as follows. 
Let $\goth A$ be a linear space
 endowed with a sequence of bilinear
operations $\ensquare n:\goth A\otimes \goth A\to \goth A,\ n\in \Z$,
and a distinguished element $\1\in \goth A$. Let  $D:\goth A\to
\goth A$ be a linear map given by $Da = a\ensquare{-2}\1$. 
Then $\goth A$ is a vertex algebra if it satisfies  
the following conditions for any $a,b,c\in \goth A$ and $m,n\in \Z$:

\begin{itemize}
\item[V1.](Locality)
$a\ensquare n b = 0$ for $n\gg0$.
\item[V2.](Identity)
$\1\ensquare n a = \delta_{n,-1}a,\ \ 
a\ensquare n \1 = \begin{cases}
0&\text{if}\quad n\ge0,\\
\frac 1{(-n-1)!} D^{-n-1} a&\text{if}\quad n<0.
\end{cases}$
\item[V3.]
(Associativity)
\vskip-25pt
\begin{align*}
\big(a\ensquare n b\big)\ensquare m c = 
& \sum_{i\ge0} (-1)^i \binom ni a\ensquare{n-i} \big(b\ensquare{m+i}c\big)\\
 -& \sum_{i\le n} (-1)^i \binom n{n-i}
b\ensquare{m+i}\big(a\ensquare{n-i} c\big).\notag
\end{align*}
\item[V4.]\vskip-10pt(Commutation)
\vskip-15pt
$$
a\ensquare m \big(b\ensquare n c\big) - 
b\ensquare n\big(a\ensquare m c\big)= 
\sum_{i\ge 0}\binom mi \big(a\ensquare i b\big)\ensquare{m+n-i}c.
$$
\end{itemize}

This is not the shortest possible list of axioms. See the references
cited above for other equivalent
definitions.

For $n,m\ge 0$, the associativity V3 is exactly the conformal Jacoby identity C4.
Among other   properties of vertex
algebras are the formulas C2 as well as the quasisymmetry identity
C3, which holds for all integer $n$. 
So vertex algebras are a special case of conformal algebras.

There is an additional axiom which is often imposed on vertex
algebras, see \cite{dl,fhl,flm,li}. A vertex algebra $\goth A$ is
called a {\it vertex operator algebra} if
\begin{itemize}
\item[V5.]
$\goth A =\bigoplus_{n\in\Z} \goth A_n$ is graded so that 
$\1\in \goth A_0$ and 
$\goth A_i \ensquare n \goth A_j \subset \goth A_{i+j-n-1}$.
\item[V6.]
There exists an element $\upsilon \in \goth A_2$ generating the Virasoro
conformal algebra $\goth{Vir} \subset \goth A$, so that the relations 
\fl(vir) hold, where $\c = c\1$ for some $c\in \Bbbk$. Also,  $\upsilon
\ensquare 0 a = Da, \ \upsilon \ensquare 1 a = (\deg a)\, a$ for all
homogeneous $a\in\goth A$.
The number $2c$ is called the {\it conformal charge} of $\goth A$. 
\end{itemize}

In \sec{series} we show that conformal and vertex algebras are
precisely the algebraic structures formed by certain formal infinite series.

\subsection{Formal series and locality}\label{sec:locality}
Let $L$ be a Lie algebra.  Denote by $L\{z\}$ the
$\Bbbk$-linear span of all series of the form 
$$
\sum_{n\in \lambda+\Z}\alpha(n)\, z^{-n-1}, \qquad \alpha(n)\in L,
\quad \lambda \in \Bbbk. 
$$  
For a linear space $V$, let $\F\{V\}\subset
(\op{gl}V)\{z\}$ be the space of all such series with coefficients in
the Lie algebra $\op{gl}V$ with the property that $\alpha(n)v = 0$ for
$n\gg0$ for any fixed $v\in V$. We call  $\F\{V\}$ the space
of {\it generalized vertex operators}. It contains the space $\F(V)$ of
ordinary vertex operators, that involve only integer powers of $z$. 
Denote by $\1\in\F(V)$ the vertex operator with the only non-zero
coefficient being $\1(-1)=\op{id}$.

The space $L\{z\}= \bigoplus_{[\lambda]\in\Bbbk/\Z}L\{z\}^{[\lambda]}$ 
is graded by the group $\Bbbk/\Z$ so that $L\{z\}^{[\lambda]}$ is the
space of all series of the form $z^{-\lambda} \alpha_\lambda(z)$, \ 
$\alpha_\lambda(z) \in L[[z^{\pm1}]]=
L\{z\}^{[0]}$. The space of vertex operators $F\{V\}$ is
a homogeneous subspace of $(\op{gl}V)\{z\}$ and  $\F(V)=\F\{V\}^{[0]}$. 

A pair of series $\alpha, \beta \in L\{z\}$ are said to be {\it local of
order} $N\in\Z_+$ \cite{dl,kac2} if 
$(z-w)^N \ad{\alpha(w)}{\beta(z)} = 0$,
or, equivalently,
$$
\sum_{s=0}^N (-1)^s \binom Ns \ad{\alpha(n-s)}{\beta(m+s)}=0
\qquad \forall\, m,n\in\Bbbk.
$$
The same applies for $\alpha,\beta \in \F\{V\}$. It is easy to see
that if series $\alpha=\sum_\lambda z^{-\lambda}\alpha_\lambda$,\  
$\beta=\sum_\lambda z^{-\lambda}\beta_\lambda\in L\{z\}$ are local of
order $N$, then any two homogeneous components $\alpha_\lambda,
\beta_\mu\in L[[z^{\pm1}]]$ are local of the same order. 

\begin{Rem}
The space $\F(V)$ of vertex operators over $V$ is a linear space over
the field $F=\Bbbk((z))$ of formal power series in $z$. Let $E \supset
F$ be the field extension of $F$ generated by $z^\lambda$ for all
$\lambda \in \Bbbk$.  Both $F$ and $E$ are
differential fields, and  $E$ is a
differential Galois extension of $F$. Then $\F\{V\}$ is a
linear space over $E$, and in fact $\F\{V\} = \F(V)\otimes_FE$. 
We observe that the multiplication by the elements
of  $E$ does not affect the locality of vertex operators. 

It is possible to consider a more general extension $E$ of $\Bbbk((z))$, for
example involving $\log z$. Most of the results in this
paper can be generalized to this more general setting.
\end{Rem}

\subsection{The products of formal series}\label{sec:products}
Recall that for two ordinary vertex operators $\alpha,\beta \in \F(V)$, one can
define products $\alpha\ensquare n \beta \in \F(V)$, \  $n\in\Z$, in
the following way.  Let $\imath_{w,z}:\Bbbk(w,z)\to \Bbbk((w\inv,z))$ and  
$\imath_{z,w}:\Bbbk(w,z)\to\Bbbk((w,z\inv))$ be the expansions of a
rational function into Laurent series at $(w,z)=(\infty,0)$ and 
$(w,z)=(0,\infty)$ respectively, so that
\begin{equation}\label{fl:imath}
\begin{split}
\imath_{w,z}(w-z)^n &=\sum_{i\ge0}(-1)^{n+i} \binom ni\, w^{n-i}z^i, \\
\imath_{z,w}(w-z)^n &=\sum_{i\ge0}(-1)^i \binom ni\, w^i z^{n-i}.
\end{split}
\end{equation}
Of course, if $n\ge 0$ then $\imath_{w,z}(w-z)^n=\imath_{z,w}(w-z)^n$.
We define
\begin{equation}\label{fl:verprod}
\big(\alpha \ensquare{n} \beta\big)(z) = 
\op{Res}_w \biggl(
\alpha(w)\beta(z)\,\imath_{w,z}(w-z)^n 
-\beta(z)\alpha(w)\,\imath_{z,w}(w-z)^n
\biggr).
\end{equation}

The $m$-th coefficient of $\alpha \ensquare{n} \beta$ is given by 
\begin{equation}\label{fl:expl}
\begin{split}
\big(\alpha \ensquare{n} \beta\big)(m) 
 = &\sum_{s\ge 0} (-1)^s \binom ns \alpha(n-s)\beta(m+s)\\
 - &\sum_{s\le n} (-1)^s \binom n{n-s} \beta(m+s)\alpha(n-s).
\end{split}
\end{equation}

If $n\ge0$ then the products $\alpha\ensquare n \beta$ make sense for
formal series $\alpha,\beta\in L[[z^{\pm1}]]$ as well. In this case
\fl(expl) simplifies to 
\begin{equation*}%\label{fl:posexpl}
\big(\alpha \ensquare{n} \beta\big)(m) 
=\sum_s (-1)^s \binom ns \ad{\alpha(n-s)}{\beta(m+s)}.
\end{equation*}
One can solve these equations with respect to
the commutators and thus recover the bracket in
$L$ from the products in $L[[z^{\pm1}]]$:
\begin{equation}\label{fl:lie}
\ad{\alpha(m)}{\beta(n)} = \sum_s \binom ms 
\big(\alpha\ensquare s \beta\big)(m+n-s)
\end{equation}
for every $m\in\Z_+,\ n\in \Z$. If $\alpha$ and $\beta$ are local,
then \fl(lie) holds for all $m,n\in \Z$.

On the other hand, the $-1$-st product is given by
$\alpha \ensquare{-1} \beta = \alpha_-\beta +
\beta\alpha_+$, where $\alpha_\pm(z) = \sum_{n\gtrless 0} \alpha(n)\,z^{-n-1}$.
It is also called {\it the normally ordered product}
and is sometimes denoted by $\:\alpha\beta\:$. For any $n<0$, we have 
$\alpha\ensquare n \beta = \frac 1{(-n-1)!} \:(D^{-n-1}\alpha)\beta\:$,
where $D=d/dz:\F\{V\}\to\F\{V\}$ is the operator of differentiation. 
It follows that vertex operators satisfy the relations V2.

Now we show how to expand these products to $\F\{V\}$ and $L\{z\}$.
For homogeneous  $\alpha = z^\lambda \alpha_0$ and 
$\beta = z^\mu \beta_0$, where $\lambda,\mu \in \Bbbk$ and $\alpha_0,
\ \beta_0$ involve only integer powers of $z$, we set 
\begin{equation}\label{fl:lambdaprod}
\alpha \ensquare n \beta = \sum_{j\ge 0} \binom \lambda j\, 
\big(\alpha_0 \ensquare{n+j} \beta_0\big)\, z^{\mu+\lambda-j}.
\end{equation}
Here 
$$
\binom \lambda j = \frac{\lambda(\lambda-1)\cdots(\lambda-j+1)}{j!}.
$$
Note that the summation in \fl(lambdaprod) is finite due to the
locality of $\alpha$ and $\beta$, and hence of $\alpha_0$ and $\beta_0$.
One can easily check that the products \fl(verprod) 
satisfy the identity \fl(lambdaprod) if we substitute  $\alpha =
z^\lambda \alpha_0$, \  
$\beta = z^\mu \beta_0$ for $\lambda,\mu\in \Z$.
We extend the definition of the products $\ensquare n$ by 
linearity to arbitrary $\alpha,\beta \in \F\{V\}$ and $\alpha, \beta
\in L\{z\}$ when $n\ge 0$. Note that 
$\alpha\ensquare n \beta \in \F\{V\}^{[\lambda+\mu]}$ if 
$\alpha\in \F\{V\}^{[\lambda]}$ and 
$\beta \in \F\{V\}^{[\mu]}$.

It is clear that if $\alpha$ and $\beta$ are local of order $N$ then 
$\alpha \ensquare n \beta = 0$ for $n\ge N$.
It could also be shown that \fl(lie) remains valid for generalized series
as well if $\alpha$ and $\beta$ are homogeneous and  
$m\equiv\deg\alpha\mod \Z$,\ \ $n\equiv\deg \beta \mod \Z$. 

We now write explicitly the formula for the products of twisted
vertex operators $\alpha, \beta \in \F\{V\}$, analogous to \fl(verprod). 
Let $N$ be the order of locality of $\alpha$ and $\beta$. Denote 
$$
\Delta(w,z) = \sum_\lambda \sum_{j=0}^{N-n-1}\binom{-\lambda}j\,
 w^\lambda z^{-\lambda-j}\, (w-z)^j,
$$
where $\lambda$ runs over the set of degrees that
appear in $\alpha$.

Then 
\begin{align}\notag
\big(\alpha \ensquare n \beta \big)(z) 
&= \sum_\lambda \big(z^{-\lambda} \alpha_\lambda\big) \ensquare n \beta
\\ \label{fl:lprod0}
& = \sum_\lambda\sum_{j\ge 0}  \binom{-\lambda} j\,  
\big(\alpha \ensquare{n+j} \beta\big)\, z^{-\lambda-j}
\\ 
&\hskip-50pt = \op{Res}_w \Bigg(\Big(\alpha(w)\beta(z)\,\imath_{w,z}(w-z)^n 
 -\beta(z)\alpha(w)\, \imath_{z,w}(w-z)^n\Big)\, \Delta(w,z)\Bigg)
\notag
\end{align}

It is easy to see that the identities V2 hold for every $a\in  \F\{V\}$ 
and $n\in \Z$.

\subsection{The case of finite grading}\label{sec:formula}
Of a particular interest is the case when 
all the degrees of $\alpha$ are rational numbers with common
denominator $p$, i.e. 
$$
\alpha(z) = \sum_{q=0}^{p-1} z^{-q/p}\,\alpha_q(z),\qquad \alpha_q \in \F(V).
$$
In this case we can rewrite \fl(lprod0) in a different way. 
Some calculations show that 
$$
\Delta(w,z)=\sum_{q=0}^{p-1} \sum_{j= 0}^m \binom{-q/p}j 
w^{\frac qp}z^{-\frac qp-j}\,(w-z)^j
= \(\frac{w-z}{w^{\frac 1p} - z^{\frac 1p}}\)^{m+1} F_p(m),
$$
where 
\begin{multline*}
F_p(m) = \sum_{l=1-p}^{1-p+m}\( \sum_{q=0}^{p-1} \ \sum_{k\ge 0} 
(-1)^{l+q+kp}\binom{-q/p +k}{m}\binom{m+1}{l+q-kp}\)\\
\times w^{\frac{m-l+1}p -1} z^{\frac lp -m}.
\end{multline*}
Using this we get
\begin{multline}\label{fl:lprod}
%\begin{split}
\big(\alpha \ensquare n \beta \big)(z) 
= \op{Res}_w \Bigg( \Big(\alpha(w)\beta(z)\,\imath_{w,z}(w-z)^n - 
\beta(z)\alpha(w)\, \imath_{z,w}(w-z)^n\Big)\\
\times 
\(\frac{w-z}{w^{\frac 1p} - z^{\frac 1p}}\)^{N-n} F_p(N-n-1)\Bigg).
%\end{split}
\end{multline}

We also remark that the polynomial $F_p(m)$ has the following
property:
\begin{equation}\label{fl:F}
F_p(m)\big|_{w^{\frac 1p}=z^{\frac 1p}} = p^{-m}\,z^\frac{(m+1)(1-p)}p.
\end{equation}

\subsection{Algebras of formal series}\label{sec:series}
It is well known \cite{kac2,li} that a subspace $\goth A
\subset \F(V)$ of pairwise local vertex operators such that $\1\in\goth A$ and
$\goth A\ensquare n \goth A\subset \goth A$ for all $n\in\Z$ is a
vertex algebra. Similarly, a subspace $\goth L \subset L[[z^{\pm1}]]$
of local formal series such that $\goth L \ensquare n \goth L \subset
\goth L$ and $DL \subseteq L$ is a conformal 
algebra. Moreover, all vertex and conformal algebras are obtained in
that way. Now we generalize this result to the case of generalized vertex operators
and series.

\begin{Thm}\label{thm:vertw}\sl
\begin{enumerate}
\item\label{series:ver}
Let $\goth A \subset \F\{V\}$ 
be a subspace of pairwise local generalized          
vertex operators such that $\1\in \goth A$ and $\goth A \ensquare n \goth A
\subseteq \goth A$ for all $n\in \Z$. Then $\goth A$ is a vertex algebra.
\item\label{series:conf}
Let $\goth L \subset L\{z\}$
be a subspace of pairwise local generalized
series with coefficients in a Lie algebra $L$, such that $D\goth L
\subseteq \goth L$ and $\goth L \ensquare n \goth L
\subseteq \goth L$ for all $n\in \Z$. Then $\goth L$ is a conformal algebra.
\end{enumerate}
\end{Thm}

\begin{proof} We will prove \ref{series:ver}, the statement
\ref{series:conf} is proved in the same way. 
We have to show  that $\goth A$ satisfies the identities V3 and V4,
since V1 and V2 hold by assumption. 
These identities are linear combinations of ``vertex monomials'' 
of the form $a_1\ensquare{n_1}\ldots\ensquare{n_{l-1}}a_l$ (with some
order of parentheses) where  $a_i$'s are either equal to $\1$ or are
formal variables. We have to show that for any specification
$a_i\in\goth A$ the identity vanishes. 
Note that these identities are multilinear.

Let $\goth A^{\text{gr}} \subset \F\{V\} $ be the
graded closure of $\goth A$, i.e. the minimal graded subspace of
$\F\{V\}$ containing $\goth A$. Since all  homogeneous components
of the vertex operators from $\goth A$ are pairwise local, the space  
$\goth A^{\text{gr}}$ 
satisfies all the assumptions of the theorem, so we can assume that
$\goth A=\goth A^{\text{gr}}$ is graded.

Consider all the vertex operators $\beta\in  \F(V)$ such that 
$z^\lambda \beta \in \goth A$ for
some $\lambda \in \Bbbk$. Every two such vertex operators are local, 
so since the theorem is known to be true in the case of ordinary
vertex operators, these $\beta$'s
generate a vertex algebra  $\goth B \subset \F(V)$.

Let $R(a_1, \ldots, a_l)=0$ be an identity which we have to check. 
Since $R$ is multilinear, 
it is enough to check it for homogeneous $a_i = z^{-\lambda_i}b_i$, where
$b_i \in \goth B, \ \lambda_i \in \Bbbk,\ 1\le i\le l$. Sometimes we
must set $b_i=\1$ and $\lambda_i=0$. When we
substitute these expressions for $a_i$ into  $R(a_1, \ldots, a_k)=0$
and apply  \fl(lambdaprod),
we get a linear combination of vertex monomials of the form 
$$
P(\lambda_1,\ldots,\lambda_l)\, z^{-\lambda_1-\ldots-\lambda_l} z^m
\big(b_{\kappa(1)}\ensquare{n_1}\ldots\ensquare{n_{l-1}}b_{\kappa(l)}\big)
$$
where $P(\lambda_1,\ldots,\lambda_l)$ is a polynomial in 
$\lambda_1,\ldots,\lambda_l,\ m\in\Z,\ \kappa$ is a permutation of 
$\left\{1,2,\ldots l\right\}$ and the products $\ensquare{n_i}$'s are applied
according to some order of parentheses. We can cancel the common factor
$z^{-\lambda_1-\ldots-\lambda_l}$.
Now we observe that the only remaining factor in these
monomials that depends on $\lambda_i$'s is $P$, therefore, for fixed 
$b_1, \ldots, b_l$ the map $(\lambda_1, \ldots, \lambda_l)\mapsto 
R(z^{-\lambda_1}b_1,\ldots, z^{-\lambda_l}b_l)$ is a polynomial map
from $\Bbbk^l$ to $\F\{V\}$. But when all $\lambda_i\in \Z$ this map is
equal to 0, because then $a_i\in \F(V)$ generate a vertex algebra
which satisfies the identity $R$, therefore, since $|\Bbbk|=\infty$, 
$R$ must be identically 0. 
\end{proof}

\begin{Rem}
There is a more conceptual argument illustrating \thm{vertw}. 
Recall that $\F(V)$ is a linear space over the field $F = \Bbbk((z))$,
and $\F\{V\}$ is a linear space over the field $E$ 
generated over $F$ by $z^\lambda$'s, see Remark at the end of \sec{locality}.
Let $\goth B \subset \F(V)$ be the vertex algebra generated over $F$ by all 
vertex operators $\beta$ such that $\beta z^\lambda$
belongs to the graded closure of $\goth A$ for some $\lambda \in
\Bbbk$. Recall that a commutative algebra $A$ with derivation $d$ can
be thought of as a 
vertex algebra if we set $a\ensquare n b = 0$ for $n\ge 0$ and 
$a\ensquare n b = \frac 1{(-n-1)!}\,d^{-n-1}(a)b$ for $n<0$, so we can
treat both $E$ and $F$ as vertex algebras. Recall also that there is a
notion of tensor product of vertex algebras \cite{bor,kac2}, which can easily be
generalized for the case of vertex algebras over some ground vertex
algebra, so that one can consider the vertex algebra $\goth B\otimes_F E$. 
The idea is that $\goth A$ is a subalgebra of $\goth B\otimes_F E$.

Note that this technique works for a more general field extension $E$
of $\Bbbk((z))$.
\end{Rem}

\begin{Rem}
The products of the generalized vertex operators were first introduced by Li in
\cite{litwisted}. However, he deals only with generating functions of
these products, formulas \fl(lprod0) and \fl(lprod) seem to be new.
Li also proves \thm{vertw}(a) using more straightforward techniques. 
Realizations of conformal algebras by generalized series were
mentioned by Kac in \cite{kac_fd}.
\end{Rem}

\subsection{Differentiation of vertex operators}
\label{sec:D}
Assume that we are given a pair $(V, D)$ consisting of a linear
space $V$ and a linear map $D:V\to V$. Let $\Delta:\F\{V\}\to \F\{V\}$
be the linear operator on the space of generalized vertex operators 
defined by $\Delta \phi = \frac d{dz}\phi - \ad D\phi$.
A vertex operator $\phi \in \F\{V\}$ is said to be of weight $\lambda\in\Bbbk$
if $\Delta\phi = \lambda\, z\inv \phi$. Denote by $\F\{V\}_\lambda$ the
space of all vertex operators of weight $\lambda$. It is easy to see that if
$\phi\in\F\{V\}_\lambda$ then $z^\mu \phi \in \F\{V\}_{\lambda+\mu}$.

In general it is not true that any vertex operator can be represented as a sum
of homogeneous vertex operators. However the following is true.

\begin{Prop}\label{prop:delta}\sl
Let  $S\subset \F\{V\}$ be a set of pairwise local vertex operators, and let
$\goth A\subset \F\{V\}$ be the vertex algebra generated by $S$. If
$S\subset  \F\{V\}_0$ then also $\goth A\subset  \F\{V\}_0$ and if 
$S\subset \bigoplus_{\lambda\in\Bbbk}\F\{V\}_\lambda$ then also 
$\goth A\subset\bigoplus_{\lambda\in\Bbbk} \F\{V\}_\lambda$.
\end{Prop}

\begin{proof}
First of all we note that $\Delta$ is a derivation of all products: 
$$
\Delta\big(\alpha\ensquare n \beta\big) = (\Delta\alpha) \ensquare n
\beta + \alpha \ensquare n (\Delta\beta),
$$
because so are both 
$\op{ad} D$ and $\frac d{dz}$. Therefore if $\Delta\alpha =
\Delta\beta = 0$ then $\Delta\big(\alpha\ensquare n \beta\big)=0$
and the first statement follows. For the second statement, it is
enough to assume that all the generators from $S$ are homogeneous. A
pair $\alpha, \beta \in S$ can be written as $\alpha = z^\lambda
\alpha_0, \ \beta = z^\mu \beta_0$ for some $\lambda,\mu\in\Bbbk$
and $\alpha_0, \beta_0 \in \F\{V\}_0$. Now the statement follows
from the formula \fl(lambdaprod).   
\end{proof}

\subsection{Modules over vertex algebras}
\label{sec:modules}
Now we give several definitions of modules for vertex algebras.
We will call a vector space $V$ a {\it module} over a vertex algebra
$\goth A$ if there is a vertex algebra homomorphism 
$\pi:\goth A\to \F(V)$. In other words, for any $a,b\in\goth A$, the
vertex operators $\pi(a)$ and $\pi(b)$ are local and $\pi(a\ensquare n b) =
\pi(a)\ensquare n \pi(b)$, \ $\pi(\1)=\1$. 
We remark that sometimes what we call a module is called a weak module.

If instead of a homomorphism $\pi:\goth A\to \F(V)$ we have
a homomorphism $\pi:\goth A\to \F\{V\}$ of $\goth A$ into the space of
generalized vertex operators over $V$, then $V$ is called a {\it generalized module}.

Assume now that a vertex algebra $\goth A = \bigoplus_{[\lambda]\in\Gamma}
\goth A^{[\lambda]}$ is graded by a group $\Gamma\subset \Bbbk/\Z$ so
that $ \goth A^{[\lambda]} \ensquare n \goth A^{[\mu]} \subset 
\goth A^{[\lambda+\mu]}$.
Then a generalized module $V$ is called {\it twisted} if the
representation homomorphism
$\pi:\goth A\to\F\{V\}$ is homogeneous, that is, 
$\pi\goth A^{[\lambda]}\subset \F\{V\}^{[\lambda]}$. 
This definition is due to Li \cite{litwisted}. Equivalently,
twisted modules can be defined using the so-called twisted Jacoby
identity, see e.g. \cite{dongtwisted,dlm,ffr}.

\begin{Rem}
We can generalize the definition of twisted representation to the case
when $\F\{V\} =F(V)\otimes_FE$ for an arbitrary differential Galois
extension $E$ of the field $F=\Bbbk((z))$, see Remark at the end of
\sec{locality}. Note that the Galois group $\op{Gal}(E/F)$ acts on
$\F\{V\}$ in a natural way. 
Let $\Gamma \subset \op{Aut} \goth A$ be a group of
automorphisms of the vertex algebra $\goth A$ and fix a group
homomorphism $\rho:\Gamma\to \op{Gal}(E/F)$. A representation 
$\pi :A\to \F\{V\}$ is called twisted if it is equivariant with the
action of $\Gamma$: \ $\pi(\gamma a) = \rho(\gamma)\pi(a)$ for any
$\gamma \in \Gamma$ and $a\in \goth A$. Note that for different
homomorphisms $\rho:\Gamma\to \op{Gal}(E/F)$ we will get different
categories of twisted representations.
\end{Rem}

Assume that an $\goth A$-module $V$ (generalized or not) has a linear
map $D:V\to V$. Then the module $V$ is called {\it strong} if
$\pi(\goth A)\subset \F\{V\}_0$, see \sec{D}.
We note that if $\goth A$ contains an element $\upsilon$ such that 
$Y(\upsilon)(0) =D$ (as it is the case when $\goth A$ is a vertex
operator algebra), and $V$ is a module over $\goth A$ such that
$\pi(\upsilon)\in \F(V)=\F\{V\}^{[0]}$, then $V$ is a strong module. Indeed,  by
\fl(lie) we have $\ad{\upsilon(0)}{a(n)} = \big(\upsilon \ensquare 0
a\big)(n) = \big(Da\big)(n)$ for any $a\in\goth A^{[n]}$.

Let again $V$ be a module (generalized or not) over $\goth A$ and let
$\pi:\goth A\to \F\{V\}$ be the representation map. Assume
$\goth A$ is a vertex operator algebra and let $\upsilon \in \goth A$
be the Virasoro element. We say that $V$ is a module
over the vertex {\it operator} algebra if 
$V= \bigoplus_{n\in\Bbbk}V_n$ is graded, 
$\pi(\upsilon)\in \F(V)$ and 
$\pi(\upsilon)(1)\big|\raisebox{-3pt}{$V_n$} = n$. As it was mentioned
above, in this case $V$ is necessarily a strong module over $\goth A$.

\begin{Rem}
Let $\goth A$ be a vertex algebra and let $V$ be a strong twisted module over
$\goth A$. Then the semidirect product $\goth A\ltimes V$ has a
structure of {\it generalized vertex algebra}, introduced by Dong and
Lepowsky \cite{dl}. The products in a generalized vertex algebra are
indexed not necessarily by integers. If $a\in \goth A$ and 
$v\in V$ then $a\ensquare n v = a(n)v$ and the products 
$v\ensquare n a$ are defined using the quasisymmetry identity C3, see
\sec{def}. 
\end{Rem}

%\vfill
\section{Twisted realizations of conformal algebras}\label{sec:2}
In this section we consider realizations of conformal algebras by
generalized formal series. We construct the twisted coefficient
algebra of a conformal algebra, which gives the universal 
realization of this type. As an example we consider affine conformal
algebras in \sec{affine}.  
\subsection{The coefficient algebra}
\label{sec:coeff}
Let $\goth L$ be a conformal algebra. Assume that 
$\goth L=\bigoplus_{[\lambda]\in\Gamma} \goth L^{[\lambda]}$ is graded
by a group $\Gamma \subset\Bbbk/\Z$ so that 
$\goth L^{[\lambda]}\ensquare n\goth L^{[\mu]}\subseteq \goth L^{[\lambda+\mu]}$,
\ $D\goth L^{[\lambda]}\subseteq\goth L^{[\lambda]}$.

% We remark that 
% if $\chi:\Gamma\to \Bbbk^\times$ is a character of $\Gamma$,
% then one can define a conformal algebra automorphism 
% $\sigma_\chi:\goth L\to \goth L$
% such that $\sigma_\chi\big|\raisebox{-4pt}{$\goth L^{[\lambda]}$} =
% \chi(\lambda)$. If $\Gamma\cong \Z/p\Z$ then $\sigma_\chi^p=1$. 
% If $\Bbbk=\C$ then there is always a character of
% $\Gamma\subset \Bbbk/\Z$ given by  $\lambda\mapsto \exp(2\pi i\, \lambda)$. 

Define a Lie algebra $\cff_\Gamma\goth L$ in the following
way. The underlying linear space of $\cff_\Gamma\goth L$ is spanned
by the symbols $a(n)$ for all homogeneous $a \in \goth L$ and 
$\Bbbk\ni n\equiv \deg a \mod \Z$
subject to the linear relations 
$(Da)(n)=-n\,a(n-1)$. The brackets in $\cff_\Gamma\goth L$ are defined
by the formula \fl(lie) for $a=\alpha\in\goth L^{[m]}$ and 
$b=\beta\in \goth L^{[n]}$.
%\begin{equation*}%\label{fl:coeff}
%\ad{a(m)}{b(n)} = \sum_{s\ge0}\binom ms\big(a\ensquare s b\big)(m+n-s),
%\end{equation*}
For a non-homogeneous $a\in\goth L$  denote $a(n) = a^{[n]}(n)$,
where $a^{[n]}$ is the projection of $a$ onto the space $\goth
L^{[n]}$.

This construction generalizes the construction of usual coefficient
algebra $\cff\goth L = \cff_{[0]}\goth L$ done in
\cite{kac2,kac_fd,freecv}. If $\goth L$ is a vertex algebra, then
$\cff_\Gamma\goth L$ was considered in \cite{dlm}, at least when
$\Gamma$ is a finite cyclic group. 

Let $\phi\mspace{-1mu}:\goth L \to (\cff_\Gamma\goth L)\{z\}$ be a map given by
$\phi(a) = \sum_{n\in\Gamma+\Z}\mspace{-1mu} a(n)\mspace{1.5mu} z^{-n-1}$. 
It is easy to see
that $\phi$ is a homomorphism of $\goth L$ into a conformal subalgebra
of $(\cff_\Gamma\goth L)\{z\}$ such that  
$\phi\big(\goth L^{[\lambda]}\big)\subset(\cff_\Gamma\goth L)\{z\}^{[\lambda]}$. 
Moreover, this map is universal in the following sense: if 
$\rho:\goth L \to K\{z\}$ is another homogeneous homomorphism of $\goth L$
into a conformal subalgebra of series with coefficients in some Lie
algebra $K$, then there is a unique Lie algebra homomorphism 
$\pi:\cff_\Gamma\goth L\to K$ making the following diagram
commutative:
\begin{equation}\label{fl:diagram}
\begin{array}{c}
(\cff_\Gamma\goth L)\{z\} \xrightarrow{\hskip5pt\pi\hskip5pt} K\{z\} \\[3pt]
\hskip25pt 
\text{\Large $\nwarrow$}\mspace{-6mu}\raisebox{5pt}{$\scriptstyle{\phi}$} 
\hskip15pt
\text{\Large $\nearrow$}\mspace{-25mu}\raisebox{5pt}{$\scriptstyle{\rho}$}\\
\hskip33pt \goth L 
\end{array}
\end{equation}

\vskip-5pt
The proof of the fact that $\cff_\Gamma\goth L$
is indeed a Lie algebra and of the above universality property is done
in the same way as in the case when $\Gamma=0$.

We also give another construction of $\cff_\Gamma\goth L$. Let
$L=\cff\goth L$ be the ordinary coefficient Lie algebra of $\goth L$.
Let $\goth L_\Gamma \subset L\{z\}$ be the conformal algebra generated
by all the series of the form $z^n a$ for $n\in\Bbbk$ and $a\in\goth L^{[n]}$. 
The algebra $\goth L_\Gamma$ is closed under the
multiplication by $z^{\pm1}$, in fact we have 
$z^{-\lambda}\goth L_\Gamma^{[\lambda]} = 
\Bbbk[z^{\pm1}] \goth L^{[\lambda]}$. 

Consider the coefficient Lie algebra 
$\cff\goth L_\Gamma$ of $\goth L_\Gamma$.
Let $L_\Gamma\subset \cff\goth L_\Gamma$ be its subalgebra
consisting of all
elements of the form $a(0)$ for $a\in \goth L_\Gamma$. In other
words, $L_\Gamma = \op{Ker}D$, where $D:\cff\goth L_\Gamma\to
\cff\goth L_\Gamma$ is the derivation given by
$a(n)\mapsto -n\,a(n-1)$. There is a Lie algebra homomorphism
$\cff\goth L_\Gamma\to L_\Gamma$ given by $a(n) \mapsto (z^n a)(0)$
for $a\in \goth L^{[n]}$,
which induces the conformal algebra homomorphism 
$\eta:\goth L_\Gamma \to L_\Gamma[[w^{\pm1}]]$ such that 
%\begin{equation}\label{fl:eta}
$\eta(z^k a) = w^k\eta(a)$
%\end{equation}
for $k\in \Z$. It could be shown that $\eta$ is the universal
among all realizations
of $\goth L_\Gamma$ by integral formal series that commute with
multiplication by the variable.

We construct a homomorphism $\rho:\goth L \to L_\Gamma\{w\}$ in the
following way. Let $a \in \goth L^{[\lambda]}$, define
$\rho(a) = \eta(z^\lambda a)w^{-\lambda}\in L_\Gamma\{w\}$.
Clearly, $\rho$ is a linear map and it does not depend on the choice of
representative $\lambda\in[\lambda]$. Now we show that $\rho$ is
indeed a homomorphism. 

Let $a_1, \ldots, a_k\in\goth L$ be homogeneous elements of $\goth
L$ satisfying a  conformal identity 
$R(a_1, \ldots, a_k) =0$. Here $R$ is a linear combination of
conformal monomials in $a_i$'s with coefficients in
$\Bbbk[D]$. We have to show that $R\big(\rho(a_1),\ldots, \rho(a_k)\big)=0$ in 
$L_\Gamma\{w\}$. Without loss of generality we can assume that $R$ is
homogeneous with respect to the gradation by $\Gamma$.
 Let $\deg a_i=[\lambda_i]$.
Substitute $a_i = z^{-\lambda_i} b_i$ in $R$, apply the
formula \fl(lambdaprod) and then cancel the common factor 
$z^{-\lambda_1-\ldots-\lambda_k}$. We get an identity $R_1(b_1, \ldots, b_k)=0$,
which holds in $\goth L_\Gamma\subset L\{z\}$,
where $R_1$ is a combination of conformal monomials in
$b_1,\ldots,b_k$ with coefficients in $\Bbbk[D,z^{\pm1}]$. Since
$\eta: \goth L_\Gamma \to L_\Gamma[[w^{\pm1}]]$ is a
$\Bbbk[D,z^{\pm1}]$-module homomorphism, we have $R_1\big(\eta(b_1),
\ldots, \eta(b_k)\big)=0$ in $L_\Gamma[[w^{\pm1}]]$. Substitute now
$\eta(b_i) = w^{\lambda_i} \rho(a_i)$ and apply \fl(lambdaprod)
again. After dividing by $w^{\lambda_1+\ldots+\lambda_k}$,
we get $R\big(\rho(a_1),\ldots, \rho(a_k)\big)=0$. 

By the universality of $\cff_\Gamma\goth L$ we get a map 
$\pi:\cff_\Gamma\goth L\to L_\Gamma$ such that the diagram \fl(diagram)
commutes for $K= L_\Gamma$.

\begin{Prop}\label{prop:coeff}\sl
The map $\pi:\cff_\Gamma\goth L\to L_\Gamma$ is an isomorphism.
\end{Prop}
\begin{proof}
We will show that the homomorphism $\rho:\goth L\to L_\Gamma\{w\}$
constructed above has the same universality property as
$\cff_\Gamma\goth L$. Then the fact that 
$\pi :\cff_\Gamma\goth L\to L_\Gamma$ is an isomorphism follows from
uniqueness of $\cff_\Gamma\goth L$.

Let $\psi:\goth L \to K\{w\}$ be a homogeneous homomorphism of 
$\goth L$ into a space of generalized series with coefficients in some
Lie algebra $K$. This induces a homomorphism 
$\theta:\goth L_\Gamma \to K[[w^{\pm1}]]$ defined by 
$\theta(z^n a) = w^n\psi(a)$ for any $a\in\goth L^{[\lambda]}$ and
$n\in\lambda+\Z$. It is easy to see that $\theta(z^n a) = w^n
\theta(a)$ for every $a\in\goth L_\Gamma$ and $n\in\Z$. Hence there is
a homomorphism $L_\Gamma \to K$ such that all the corresponding
diagrams commute.
\end{proof}
\subsection{Twisted affine algebras}\label{sec:affine}
We illustrate the construction of \sec{coeff} by the example of
affine conformal algebras. 

Let $\goth g$ be a  Lie algebra  with an
invariant bilinear form $(\,\cdot\ |\ \cdot\,)$, i.e. such
that $(\ad ab\,|\,c)=(a\,|\,\ad bc)$. Assume that the algebra 
$\goth g = \bigoplus_{[\lambda]\in\Gamma} \goth g^{[\lambda]}$  is graded
by a group $\Gamma \subset\Bbbk/\Z$. Assume further that
the gradation on $\goth g$ agrees with the form   
$(\,\cdot\,|\,\cdot\,)$ in the following way: 
\begin{equation}\label{fl:grform}
\big(\goth g^{[\lambda]}\big|\goth g^{[\mu]}\big)=0 \quad \text{unless}\quad
\lambda+\mu=0.
\end{equation}
Consider the Lie algebra $\^L_\Gamma=\goth g \otimes \Bbbk[\Gamma+\Z]
\oplus \Bbbk \c$, where $\Bbbk[\Gamma+\Z]$ is the group algebra of 
$\Gamma+\Z$. We will write $a(n) = a\otimes n$ for $a\in \goth g, \ 
n \in \Gamma+\Z$. The brackets in $\^L_\Gamma$ are defined by
$$
\ad{a(m)}{b(n)} = \ad ab (m+n) +
\delta_{m+n,0}\ m\, \big(a^{[m]}\big|b^{[n]}\big)\, \c, \qquad
\c \in Z\big(\^L_\Gamma\big),
$$
where $a^{[m]}$ and $b^{[n]}$ are projections of $a$ and $b$ onto 
$\goth g^{[m]}$ and $\goth g^{[n]}$ respectively.
The twisted affine Lie algebra 
$L_\Gamma\subset \^L_\Gamma$ is the subalgebra of 
$\^L_\Gamma$ spanned by $\c$ and all elements of the form $a(n)$ for
$a\in \goth g^{[n]}$.
The grading on $\goth g$ induces a grading on $L_\Gamma$ by setting 
$\deg a(n) = \deg a$, \ $\deg \c = 0$. From now on, if $a\otimes n \in
\^L_\Gamma\ssm L_\Gamma$, then we set $a(n)=0$.

For any $a\in\goth g$ consider series
$\~a = \sum_{n\in\Gamma+\Z}a(n)\,z^{-n-1}\in L_\Gamma\{z\}$. These
series are local of order 2, and together with the series $\c =
\c\,z^0$ they generate a conformal algebra $\goth L \subset
L_\Gamma\{z\}$. The products between $\~a$'s are 
\begin{equation}\label{fl:affprod}
\~a \ensquare 0 \~b = \~{\ad ab},\qquad 
\~a \ensquare 1 \~b = (a|b)\,\c.
\end{equation}

The affine conformal algebra $\goth L$ is a homogeneous subalgebra of
$L_\Gamma\{z\}$ so that for $a\in \goth g^{[\lambda]}$ we have 
$\~a \in L_\Gamma\{z\}^{[\lambda]}$. It is independent on the
$\Gamma$-grading of $\goth g$.  It is easy to see that 
$L_\Gamma= \cff_\Gamma\goth L$.

For each equivalence class $[\lambda] \in \Gamma\subset\Bbbk/\Z$ choose
a representative $\lambda\in \Bbbk$ such that the representative of
$\Z$ is 0 and if $[\lambda] + [\mu]=0$, then either $\lambda=\mu=0$ or 
$\lambda+\mu=1$. 
For example, if $\Bbbk\subseteq \mathbb R$, then we can take 
$0 \le \lambda<1$. 
A very important special case is when $\Gamma$ is a finite cyclic
group of order $p$, then one can take the set of representatives of
$\Gamma$ in $\Bbbk$ to be 
$\big\lbrace 0,\frac 1p, \frac 2p,\ldots,\frac{p-1}p\big\rbrace$.

For a homogeneous element $a\in \goth g^{[\lambda]}$ consider the
series $\tau_a = z^\lambda \~a \in L_\Gamma[[z^{\pm1}]]$.  Clearly,
the series $\tau_a$ are pairwise local of order 2, so together with
$\c$ they generate a conformal algebra $\goth L_\Gamma\subset
L_\Gamma[[z^{\pm1}]]$, called a twisted affine conformal algebra.  
We calculate the non-zero products of these
series, using \fl(lambdaprod) and \fl(affprod).  Here $a,b\in \goth g,
\ \deg a=[\lambda],\ \deg b = [\mu]$:
\begin{align*}
\tau_a \ensquare 0 \tau_b 
&=z^{\lambda+\mu}\,\~{\ad ab} + \lambda\, z^{\lambda+\mu-1}\, (a|b)\,\c\\
&=  \begin{cases}
z\ \tau_{\ad ab}+\lambda\, (a|b)\,\c&\text{if}\quad \lambda+\mu\ge1,\\
\tau_{\ad ab}&\text{if}\quad \lambda+\mu<1,
\end{cases}\\[5pt]
\tau_a \ensquare 1 \tau_b &= z^{\lambda+\mu}\, (a|b)\, \c = 
\begin{cases}
(a|b)\,\c & \text{if}\quad \lambda=\mu=0,\\
z\, (a|b)\,\c & \text{if}\quad \lambda+\mu=1,\\
0 &\text{otherwise}.
\end{cases}
\end{align*}
It is not difficult to see that $L_\Gamma$ is isomorphic to the
subalgebra of $\cff \goth L_\Gamma$ consisting of all elements of the
form $a(0)$, i.e., the embedding $\goth L_\Gamma\to
L_\Gamma[[z^{\pm1}]]$ is the universal one among realizations of
$\goth L_\Gamma$ by formal series that agree with the multiplication
by $z$, see \sec{coeff}.

\begin{Rem}
While in general we cannot guarantee that the twisted representation
map $\rho:\goth L\to L_\Gamma\{w\}$ is injective, in the particular
case when $\goth L$ is an affine conformal algebra we do know that
$\rho$ is an isomorphism.
\end{Rem}

\begin{Rem}
In principal, one can use the construction of \sec{coeff} to get 
twisted representations of a vertex algebra. Namely, let
$\goth A = \bigoplus_{[\lambda]\in\Gamma} \goth A^{[\lambda]}$ be a
graded vertex algebra, and let $\pi:\goth A \to \F(V)$ be its
representation. Let $\goth A_\Gamma \subset \F\{V\}$ be the vertex
algebra generated by the generalized vertex operators  $z^n \pi(a)$ for 
$a \in \goth A^{[n]}$. Then $\goth A_\Gamma$ will be closed under the
multiplication by $z^{\pm1}$. If we have a representation 
$\eta:\goth A_\Gamma \to \F(U)$ such that 
$\eta(z^ka) = w^k\eta(a)$ for $k\in \Z$, then the map 
$a\mapsto \eta(z^\lambda a)w^{-\lambda}$ for $a\in\goth A^{[\lambda]}$ 
defines a twisted
representation of $\goth A$ on $U$. However, in contrast with the
conformal case, it is not clear how to construct such a representation
$\eta$ of $\goth A_\Gamma$. 
On the other hand, if $\goth A$ is generated by a conformal algebra
$\goth L\subset \goth A$, then applying the construction of
\sec{coeff} to $\goth L$ we can get a twisted representation of
possibly some other enveloping vertex algebra $\goth B\supset \goth L$. 
\end{Rem}

\section{Lattice vertex algebras}\label{sec:3}
In this section apply the technique developed in \sec{1} to lattice
vertex algebras. We assume here that the ground field is $\C$ and the
group $\Gamma$ is the cyclic group of order $p\ge0$. 
We identify $\Gamma\subset\C/\Z$ with the set 
$\big\lbrace 0,\frac 1p, \frac
2p,\ldots,\frac{p-1}p\big\rbrace$.
 
We remark that most of the constructions below can be done in a much more
general setting. 
\subsection{Representation theory of Heisenberg algebras}
\label{sec:heisen}
Let $\goth h$ be a $\Bbbk$-vector space of dimension $l<\infty$
equipped with a non-degenerate bilinear form $(\,\cdot\,|\,\cdot\,)$. Assume that
$\goth h = \bigoplus_{\lambda\in \Gamma} \goth h^{[\lambda]}$ is
graded by $\Gamma$ and that the grading agrees with the form in sense of
\fl(grform).  Recall that the grading induces an automorphism
$\sigma:\goth h\to\goth h$ such that
$\sigma\big|\raisebox{-3pt}{$\goth h^{[\lambda]}$} = \exp(2\pi
i\,\lambda)$. In fact the existence of an automorphism 
$\sigma:\goth h\to\goth h$ of order $p$ is equivalent to the existence
of the above grading of $\goth h$ by the group $\Gamma = \Z/p\Z$.  
Note that $\sigma$ preserves the norm on $\goth h$. For 
$h\in \goth h$ we denote by $h^{[\lambda]}$ the projection of $h$ onto 
$\goth h^{[\lambda]}$.

View $\goth h$ as an Abelian Lie algebra, and 
let 
$$
H=H_\Gamma=\op{Span}\bigset{a(n)}{a\in \goth h^{[n]}}\oplus \C\c
$$ 
be the corresponding (twisted,  unless $p=1$) affine Lie algebra, see
\sec{affine}. It
is usually called a Heisenberg Lie algebra. As in
\sec{affine}, for an element $a\in \goth h$ consider formal series 
$\~a = \sum_{n \in\Gamma + \Z} a(n)\,z^{-n-1} \in H\{z\}$. 
These series, together with $\c$,
span over $\C[D]$ a copy of conformal Heisenberg algebra 
$\goth H \subset  H\{z\}$, so that $H = \cff_\Gamma\goth H$.
The grading on $\goth h$ lifts to a grading on $H$ and $\goth H$ and  the
automorphism $\sigma$ lifts to automorphisms of $H$ and $\goth H$. 
Recall also that there is another grading on $H$ given by setting 
$\deg a(n) = -n\in \frac 1p\Z$.

We note that $\goth h^{[0]}\subset Z(H)$ so that $H= H'\oplus \goth
h^{[0]}$, where 
$$
H' = \op{Span}\set{a(n)}{a\in\goth h, \ n\in\C^\times}\oplus \C\c.
$$ 
Let $H_{\pm}=\op{Span}\set{a(n)}{n\gtrless0}\subset H'$. 
We have $H'=H_-\oplus\C\c\oplus H_+$.

Now let $M$ be a restricted $H$-module, i.e. such that for any
$u\in M$ we have $h(n)u=0$ for $n\gg0$. Assume that $\c$ acts on $M$
by the identity. Then the vertex operators  $\~h \in \F\{M\}$ 
generate a vertex algebra $\goth V_0 \subset \F\{M\}$. 
It is an {\it enveloping vertex
algebra} of the conformal Heisenberg algebra $\goth H$. 
The algebra $\goth V_0$ is a module over the Heisenberg algebra $H$ by 
$h(n)x = \~h\ensquare n x$.  It is well known that $\goth V_0$ is the
unique 
enveloping vertex algebra of $\goth H$ such that $\c=\1$. As a module
over $H$, the vertex algebra $\goth V_0$ is isomorphic to the
so-called canonical relations representation $M(1)=U(H_-)\1$, which
is generated by a single element $\1$ such that 
$H_+\1 = 0$. 
  
The vertex algebra $\goth V_0$ is in fact
a vertex operator algebra: it is graded so that $\deg \~h = 1$ for
$h\in \goth h$ and it contains a Virasoro element 
$\upsilon=\frac12\sum_{i=1}^\ell\~\alpha_i\ensquare{-1}\~\beta_i$, 
where $(\alpha_1,\ldots,\alpha_\ell)$ and $(\beta_1,\ldots,\beta_\ell)$
are dual bases of $\goth h$, i.e. such that 
$(\alpha_i|\beta_j)=\delta_{ij}$. 
We have $\upsilon\ensquare 0 u = Du$ for 
all $v\in \goth V_0$, \  $\upsilon \ensquare 1 u = (\deg u) u$ for all
homogeneous $u\in \goth V_0$, \ $\upsilon \ensquare 2\upsilon =0$ 
and $\upsilon \ensquare 3 \upsilon = \frac 12 \1$.  

An $H$-module $V$ is called $\goth h^{[0]}$-diagonalizable if it can
be decomposed into a direct sum of subspaces 
$V = \bigoplus_{\xi\in(\goth h^{[0]})^*}V_\xi$,
so that for $h \in \goth h^{[0]}$ one has 
$h\big|\raisebox{-3pt}{$V_\xi$} = \xi(h)$.
Recall \cite{kac1} that an $H$-module $V$ belongs to the category $\cal O$ 
if $V$ is $\goth h^{[0]}$-diagonalizable and for any $v\in V$ there is
$n\in \frac 1p\Z$ such that for any $x\in U(H)$ of 
$\deg x \ge n$ we have $xv=0$. Clearly, any module from the category
$\cal O$ is restricted. 

For the future reference we cite here a result from the
representation theory of Heisenberg algebras.

\begin{Lem}\label{lem:heisen}\sl
Let $V$ be a module over the Heisenberg Lie algebra $H$. Let 
$\Omega = \set{v\in V}{H_+v=0}\subset V$ be the vacuum subspace of $V$.
Then the following conditions are equivalent:
\begin{itemize}
\item[\hbox to17pt{\rm\hfill(i)}]
$V\cong M(1)\otimes \Omega$ and $\Omega$ is $\goth h^{[0]}$-diagonalizable;
\item[\hbox to17pt{\rm\hfill(ii)}]
$V=U(H)\Omega$ and $\Omega$ is $\goth h^{[0]}$-diagonalizable;
\item[\hbox to17pt{\rm\hfill(iii)}]
$V\in\cal O$ and $V$ is completely reducible over $H$;
\item[\hbox to17pt{\rm\hfill(iv)}]
$V\in \cal O$ and there is a grading 
$\displaystyle{\smash[b]{V=\bigoplus_{n\in \frac 1p\Z} V_n}}$
on $V$ such that $\deg a(n)=-n$ and 
$\upsilon(1)\big|\raisebox{-4pt}{$V_n$} = n$, where 
$$
\upsilon(1) = \frac 12 \sum_{i=1}^l\Bigg(
\sum_{s<0}\alpha_i(s)\beta_i(-s) +\sum_{s\ge0} \beta_i(-s)\alpha_i(s)
\Bigg)\in \overline{U(H)}
$$   
is the first coefficient of the Virasoro element 
$\upsilon\in \goth V_0$.
\end{itemize}
\end{Lem}
Clearly, the vacuum space $\Omega\subset V$ is stable under the action 
of $\goth h^{[0]}$.
The condition (iv) means that $V$ is a module over the vertex operator
algebra $\goth V_0$, see \sec{modules}. We remark that if in (iv) we
assume that the grading on $V$ is bounded from below, that is, there
is $n_0\in \frac1p\Z$ such that $V_n=0$ for $n<n_0$, then the
condition $V\in\cal O$ becomes obsolete. 

\subsection{Fock spaces}
\label{sec:fock}
Let $\Lambda \subset \goth h$ be a lattice in $\goth h$ of rank
$l$. Assume that $\Lambda$ agrees with  the $\Gamma$-grading on 
$\goth h$ in the sense that $\op{rk} \Lambda^{[\lambda]}=\dim \goth
h^{[\lambda]}$ for $[\lambda] \in \Gamma$. Consider a central extension
$$
1\longrightarrow \Phi \longrightarrow \^\Lambda \longrightarrow \Lambda
\longrightarrow 1
$$
of $\Lambda$ by the multiplicative group $\Phi\cong\C^\times=\C\ssm\{0\}$.  
Let $e:\Lambda \to \^\Lambda$ be a section
and $\e:\Lambda\times\Lambda \to Z$ be the corresponding
2-cocycle, so that $e(\alpha)e(\beta) = \e(\alpha, \beta)\,
e(\alpha+\beta)$ for $\alpha, \beta \in \Lambda$.
Denote by $R$ the quotient of the group algebra
$\C\big[\^\Lambda\big]$ obtained by the  
identification of $\Phi$ with $\C^\times\subset \C$. 

The key object of the construction is a vector space $V$, which is a
module over the associative algebra $R$ and a restricted module over
the twisted Heisenberg algebra $H$ such that 
\begin{equation}\label{fl:RH}
\ad{h(n)}{e(\alpha)} = \delta_{n,0}(h|\alpha)\,e(\alpha),\quad
\text{for}\quad h\in\goth h, \quad n\in \Gamma+\Z,
\end{equation}
and $\c=\op{id}$.
The formulas \fl(RH) mean that $H$ acts on $R$ by derivations, so we
can form a skew tensor product algebra $U(H)\,\~\otimes\, R$ such that $V$ is a
$U(H)\,\~\otimes\, R$-module.

We assume that $V= \bigoplus_{\xi\in(\goth h^{[0]})^*}V_\xi$  
is $\goth h^{[0]}$-diagonalizable. 
The bilinear form on $\Lambda$ induces a homomorphism 
$\nu:\Lambda \to \big(\goth h^{[0]}\big)^*$ by $\nu(\alpha)\beta
= (\alpha|\beta)$ for $\alpha \in \Lambda$, \ $\beta \in \Lambda^{[0]}$.
The relation $\ad{\beta(0)}{e(\alpha)} = (\alpha|\beta)e(\alpha)$ is
equivalent to the fact that 
$e(\alpha)V_\xi \subseteq V_{\xi+\nu(\alpha)}$.

Such a module $V$ is sometimes referred to as a {\it Fock space}. In
\sec{lattice} we will make $V$ to be a module over certain vertex
algebra, under some additional assumptions. 

A standard way to construct a Fock space $V$ is as follows. 
Let $\Omega=\bigoplus_{\xi\in(\goth h^{[0]})^*}\Omega_\xi $ be a
$\goth h^{[0]}$-diagonalizable module over $R$ and $\goth h^{[0]}$
such that $\ad{h}{e(\alpha)} = (h|\alpha)\,e(\alpha)$ for $h\in
\goth h^{[0]}$. In other words,
$\Omega$ is a module over a certain skew tensor product algebra 
$U\big(\goth h^{[0]}\big)\,\~\otimes\, R$. As before, we have 
$e(\alpha)\Omega_\xi \subseteq \Omega_{\xi+\nu(\alpha)}$.
Take now a restricted $H$-module $M$ such that $\goth h^{[0]}$ acts on
$M$ by 0. Let $V = M\otimes \Omega$. Then 
$V=\bigoplus_{\xi\in(\goth h^{[0]})^*}V_\xi$,\ 
$V_\xi = M\otimes \Omega_\xi$, is also a restricted
$\goth h^{[0]}$-diagonalizable 
$H$-module and also it is a module over $R$ such
that the relations \fl(RH) hold. 
%The action is given by 
%$a(n)\big(v\otimes u\big) = \big(a(n)v\big)\otimes u$,\ 
%$a(0)\big(v\otimes u\big) = v\otimes \big(a(0)u\big)$,\
%$x(v\otimes u) = v \otimes x(u)$,
%for $a\in \goth h, \ 0\neq n\in \Gamma+\Z, \ x\in R, 
%\ v\in M, \ u\in \Omega$. The central element $\c$ acts on $V$ by the
%identity. 

For example, assume that  $V$ satisfies the conditions of \lem{heisen}
as a module over $H$. Then $V=M(1)\otimes \Omega$ can be obtained by
the above construction for $M=M(1)$. The vacuum space 
$\Omega \cong \1\otimes  \Omega\subset V$ becomes a module over 
$U\big(\goth h^{[0]}\big)\,\~\otimes\, R$.

\subsection{Twisted vertex operators}
\label{sec:tvo}
Now we are ready to define the main ingredient of this construction --- the
vertex operator $X_\alpha \in \F\{V\}$ for $\alpha\in \Lambda$. 
Let  $h' = h-h^{[0]}$ for $h\in\goth h$. We  set 
\begin{equation}\label{fl:VO}
X_\alpha(z) = e(\alpha)\,E_-(\alpha,z)E_+(\alpha,z)\, 
z^{\alpha(0)}\, z^{- (\alpha'|\alpha')/2},
\end{equation}
where 
$$
E_{\pm}(\alpha,z) = 
\exp \mspace{-10mu}\sum_{n\in \Gamma+\Z, \ n\gtrless 0}\mspace{-10mu}
- \frac{\alpha(n)}n\, z^{-n}.
$$
Note that $\alpha(0) = \alpha^{[0]}$ is the projection of $\alpha$
onto $\goth h^{[0]}$, so that
$z^{\alpha(0)}\big|\raisebox{-3pt}{$V_\xi$} = 
z^{\xi(\alpha^{[0]})}$.

\begin{Rem}
In fact we need that $V$ is $\goth h^{[0]}$-diagonalizable only for
the expression $z^{\alpha(0)}$ to make sense. We can instead interpret
$z^{\alpha(0)}=\exp\big(\alpha(0)\,\log z \big)$, and then this
expression is well defined under somewhat weaker assumptions, for example
it is enough to require that $\alpha(0)$ acts locally finite
dimensionally. 
\end{Rem}

\begin{Prop}{\rm\cite{lep}}\label{prop:VO}\sl\quad 
Let $h\in \goth h, \ n\in \Gamma+\Z, \ \alpha, \beta  \in \Lambda$.
We have 
\smallskip
\begin{enumerate}
\item\label{VO:ad}
$\ad{h(n)}{X_\alpha(z)}=(\alpha|h)\,z^n\,X_\alpha(z)$;\medskip
\item\label{VO:n}
$\~h\ensquare n X_\alpha = 0$ \  if  \ $1\le n\in\Z_+$ and
$\~h\ensquare 0 X_\alpha = (\alpha|h)\, X_\alpha$; \medskip
\item\label{VO:D}
$DX_\alpha = \~\alpha \ensquare{-1} X_\alpha$;\medskip
\item\label{VO:vir}
$\upsilon \ensquare 0 X_\alpha = DX_\alpha, \ \ 
\upsilon \ensquare 1 X_\alpha = \frac 12 (\alpha|\alpha)\, X_\alpha$;
\item\label{VO:XX}
$\displaystyle{
X_\alpha(w) X_\beta(z) = \e(\alpha, \beta)\, X_{\alpha, \beta}(w,z)\  
\imath_{w,z}
\prod_{s=0}^{p-1} \bigl(w^{\frac 1p} - \omega^s z^{\frac 1p}\bigr)^{(\sigma^{-s}\alpha|\beta)},}$\newline
where $\omega = \exp \frac{2\pi\, i}p$ is the primitive $p$-th root of
unity and
\begin{multline*}
X_{\alpha, \beta}(w,z) = e(\alpha+\beta)\, 
E_-(\alpha, w) E_-(\beta, z)E_+(\alpha, w) E_+(\beta, z)\\
\times w^{\alpha(0)} z^{\beta(0)} 
w^{- (\alpha'|\alpha')/2} z^{- (\beta'|\beta')/2}.
\end{multline*}
\end{enumerate}

\end{Prop}

We have $X_{\alpha, \beta}(w,z) = X_{\beta, \alpha}(z,w)$
and $X_{\alpha, \beta}(z,z) = 
X_{\alpha+\beta}(z)\, z^{(\alpha'|\beta')}$. The notation $\imath_{w,z}$
in \ref{VO:XX} is a short for $\imath_{w^{\frac 1p},z^{\frac 1p}}$,
see \fl(imath).
\begin{proof}
\ref{VO:ad}.\quad
Let $0\neq n \in \Gamma+\Z$. Then 
$$
\biggl[ h(n), \exp\(-\frac{\alpha(m)}m\, z^{-m}\)\biggr] = 
\delta_{n,-m}\, \exp\(-\frac{\alpha(m)}m\, z^{-m}\) (\alpha|h)\, z^n,
$$
and $h(n)$ commutes with all the rest of the factors in \fl(VO). Also,
$h(0)$ commutes with all the factors in \fl(VO) except $e(\alpha)$, whose
commutators are given by \fl(RH), so \ref{VO:ad} follows.

\medskip\ref{VO:n}.\quad
It follows that $\ad{h(n)}{X_\alpha(m)} = X_\alpha(m+n)$ 
for every $m,n\in \Gamma+\Z$. Hence we have for $n\in \Z_+$
\begin{multline*}
\big(\~h \ensquare n X_\alpha\big)(m) = \sum_s (-1)^s \binom ns 
\ad{h(n-s)}{X_\alpha(m+s)}\\ 
= (\alpha|h) \sum_s (-1)^s \binom ns X_\alpha(m+n)
= \begin{cases}
0&\text{if}\quad n>0,\\
(\alpha|h)\, X_\alpha(m) &\text{if}\quad n=0.
\end{cases}
\end{multline*}

\medskip\ref{VO:D}.\quad
Using that $(\alpha'|\alpha')/2 = 
\sum_{\lambda\in \Gamma}\limits\lambda \, (\alpha^{[\lambda]}|\alpha)$, we
get
\begin{equation*}
\begin{split}
DX_\alpha(z) &= \sum_{n<0}\alpha(n)\,z^{-n-1}\, X_\alpha(z) + 
X_\alpha(z)\sum_{n>0}\alpha(n)\,z^{-n-1} \\
&+ X_\alpha(z)\, \alpha(0)\,z^{-1}
-\sum_{\lambda\in\Gamma} \lambda\, (\alpha^{[\lambda]}|\alpha)\,
z\inv\, X_\alpha(z)\\
&= \:\~\alpha X_\alpha\: - 
\sum_{\lambda\in\Gamma} \lambda\, (\alpha^{[\lambda]}|\alpha)\,
z\inv\, X_\alpha(z).
\end{split}
\end{equation*}

On the other hand, set 
$\~\alpha(z) = \sum_\lambda
z^{-\lambda}\,\~\alpha^{[\lambda]}(z)$. The locality of $\~\alpha$ and
$X_\alpha$ is 1, hence, using \fl(lambdaprod),
\begin{equation*}
\begin{split}
\~\alpha\ensquare{-1} X_\alpha &= 
\sum_\lambda \big(z^{-\lambda}\~\alpha^{[\lambda]}\big)\ensquare{-1}X_\alpha\\
&= \sum_\lambda \biggl(
z^{-\lambda}\big(\~\alpha^{[\lambda]}\ensquare{-1}X_\alpha\big)
-\lambda\, z^{-\lambda-1}\big(\~\alpha^{[\lambda]}\ensquare{0}X_\alpha\big)
\biggr)\\
&= \:\~\alpha X_\alpha\: - 
\sum_{\lambda\in\Gamma} \lambda\, (\alpha^{[\lambda]}|\alpha)\,
z\inv\, X_\alpha(z).
\end{split}
\end{equation*}

\medskip\ref{VO:vir}.\quad
Let us calculate $\upsilon \ensquare 0 X_\alpha$ by the
associativity formula V3. We have, using \ref{VO:n} and
\ref{VO:D}:
\begin{align*}
\upsilon \ensquare 0 X_\alpha &= \frac 12\sum_i \(
\sum_{s<0}\alpha_i(s)\beta_i(-s-1)+\sum_{s\ge0}\beta_i(-s-1)\alpha_i(s)
\)X_\alpha \\
&=\frac 12\sum_i \Big(\alpha_i(-1)(\beta_i|\alpha) +
\beta_i(-1)(\alpha_i|\alpha)\Big)X_\alpha \\
&= \alpha(-1)X_\alpha = DX_\alpha.
\end{align*}
The other relation is proved in the same way.

\medskip\ref{VO:XX}.\quad
Let us first calculate $S = E_+(\alpha, w) E_-(\beta, z)
E_+\inv(\alpha,w)$. Since 
\begin{equation*}
\begin{split}
\exp\(-\frac{\alpha(n)}n w^{-n}\)& E_-(\beta, z) 
\exp\(\frac{\alpha(n)}n w^{-n}\) \\
&= \exp\(\frac{w^{-n}}n\op{ad} \alpha(n)\) E_-(\beta, z)\\
&= \big(\alpha^{[n]}\big|\beta\big)\ \frac{z^n w^{-n}}n \ E_-(\beta, z),
\end{split}
\end{equation*}
we have 
\begin{equation*}
\begin{split}
S &= \exp  \sum_{0<n \in \Gamma+\Z}\! - \big(\alpha^{[n]}\big|\beta\big)\, 
\frac{z^n w^{-n}}n\\
&= \exp \!\!\!\!\sum_{\lambda \in \left\{\frac 1p,\ldots, \frac{p-1}p,1\right\}} 
\!\! \!\!- \big(\alpha^{[\lambda]}\big|\beta\big) \sum_{n\in\lambda+\Z_+} \frac{z^n w^{-n}}n.
\end{split}
\end{equation*}

Denote $y = (z/w)^{\frac 1p}$ and let 
$\lambda = \frac qp, \ 1\le q\le p$. Then 
$$
\sum_{n\in\lambda+\Z_+} \frac{z^n w^{-n}}n 
= \int \frac{(z/w)^{\lambda-1}}{1-(z/w)}\ d(z/w)
= p \int \frac{y^{q-1}}{1-y^p}\ dy.
$$
Using the elementary fraction decomposition
$$
\frac{y^r}{1-y^p} = 
\frac 1p \sum_{s=0}^{p-1} \frac{\omega^{-sr}}{1-\omega^sy}, 
$$
we get that 
$$
\sum_{n\in\lambda+\Z_+} \frac{z^n w^{-n}}n 
= -\sum_{s=0}^{p-1}\ln\big(1-\omega^sy\big)^{\omega^{-sq}},
$$
so that 
\begin{align*}
S &= \prod_{s,q=0}^{p-1} \big(1-\omega^s y\big)^{\omega^{-sq}
(\alpha^{(q/p)}|\beta)}
= \prod_{s=0}^{p-1} \big(1-\omega^s y\big)^{\sum_q \omega^{-sq}
(\alpha^{(q/p)}|\beta)}\\
&= \prod_{s=0}^{p-1} \big(1-\omega^s y\big)^{(\sigma^{-s}\alpha|\beta)}.
\end{align*}
The rest of non-trivial commutation relations between the factors in
\fl(VO) are 
$$
e(\alpha)\e(\beta) = \e(\alpha, \beta)\, e(\alpha+\beta),\qquad
w^{\alpha(0)}e(\beta) = e(\beta)w^{\alpha(0)}\,
w^{(\alpha^{[0]}|\beta)}, 
$$
so we finally get, using that 
$\alpha^{[0]} = \frac 1p \sum_{s=0}^{p-1} \limits \sigma^s \alpha$,
\begin{equation*}
\begin{split}
X_\alpha(w) X_\beta(z) &= \e(\alpha, \beta)\,
X_{\alpha, \beta}(w,z)\, w^{(\alpha^{[0]}|\beta)} S\\
&= \e(\alpha, \beta)\, X_{\alpha, \beta}(w,z)\ 
\imath_{w,z}
\prod_{s=0}^{p-1} 
\big(w^\frac1p-\omega^s z^\frac1p\big)^{(\sigma^{-s}\alpha|\beta)}.
\end{split}
\end{equation*}
\end{proof}

The statement \ref{VO:n} implies that $X_\alpha$ generates a highest
weight module over the Heisenberg Lie algebra $H$, where the action is
given as usual by $h(n) u = \~h\ensquare n u$. This module is
irreducible and is in fact a module over the Heisenberg vertex algebra
$\goth V_0 = U(H)\1$. 

The Virasoro element $\upsilon \in \goth V_0$ gives the operator
\begin{equation}\label{fl:Dvir}
D=\upsilon(0) = \sum_i\(\sum_{s<0} \alpha_i(s)\beta_i(-s-1) 
+ \sum_{s\ge0}\beta_i(-s-1)\alpha_i(s) \).
\end{equation}
By \ref{VO:vir} the vertex operators $X_\alpha$ are of weight 0 with
respect to this $D$, see \sec{D}. 

\subsection{Lattice vertex algebras and their generalized representations}
\label{sec:lattice}
In this section we construct the vertex superalgebra $\goth
V_\Lambda$, known as the lattice vertex superalgebra. 

We start with the remark 
that only minimal modification is needed to transfer all the
definitions and results of \sec{1} and \sec{2} to the
realm of superalgebras. All commutators must be interpreted as
supercommutators, the formula \fl(verprod) must change to 
\begin{multline*}
\big(\alpha \ensquare{n} \beta\big)(z) = 
\op{Res}_w \biggl(
\alpha(w)\beta(z)\,\imath_{w,z}(w-z)^n \\
-(-1)^{p(\alpha)p(\beta)}\beta(z)\alpha(w)\,\imath_{z,w}(w-z)^n
\biggr),
\end{multline*}
etc. Notably, in the definition of vertex operator superalgebra the
eigenvalues of the grading derivation $\upsilon(1)$ are allowed to be
half-integer, and the vertex algebra 
$\goth A =\bigoplus_{n\in\frac12\Z} \goth A_n$ is graded by
$\frac12\Z$ so that the even and odd parts of $\goth A$ are
respectively $\goth A\even = \bigoplus_{n\in \Z} \goth A_n, \ 
\goth A\odd = \bigoplus_{n\in \Z+\frac12} \goth A_n$.
A reader who is averse to supermathematics can assume that the lattice
$\Lambda$ is even, i.e. $(\alpha|\alpha)\in2\Z$ for all 
$\alpha\in\Lambda$, see (i) below.

Here we make the following assumptions on the data introduced in
\sec{tvo}. 
\begin{itemize}
\item[(i)]
For any $\alpha, \beta \in \Lambda$ the numbers $m_s =
(\sigma^{-s}\alpha|\beta), \ 0\le s\le p-1$, are integer, in other words
$\Lambda\subset \goth h$ is
an integer lattice and the automorphism $\sigma:\goth h\to \goth h$
preserves the dual lattice $\Lambda' = 
\set{\alpha\in \goth h}{(\alpha|\Lambda)\subset \Z}\supseteq \Lambda$.
\item[(ii)]
The cocycle $\e:\Lambda\times\Lambda \to \C^\times$ is such that the
corresponding  commutator map is 
\begin{equation}\label{fl:comm}
C(\alpha, \beta)=\e(\alpha, \beta)\e(\beta, \alpha)\inv = 
(-1)^{(\alpha|\alpha)(\beta|\beta)+p(\alpha^{[0]}|\beta^{[0]})}
\,\omega^{-\sum_{s=1}^{p-1} s\, m_s}.
\end{equation}
\end{itemize}
Note that $p\,(\alpha^{[0]}|\beta^{[0]}) = \sum_{s=0}^{p-1} m_s\in \Z$.

The cocycle $\e$ satisfying \fl(comm) can be easily constructed 
in the following way. Let $\alpha_1, \ldots, \alpha_l$ be a $\Z$-basis
of $\Lambda$. Define $\e$ first for $\alpha, \beta \in \left\{\alpha_1,
\ldots, \alpha_l\right\}$ such
that \fl(comm) holds. This is possible since
$C(\alpha,\beta)=C(\beta,\alpha)\inv$ and $C(\alpha,\alpha)=1$. 
Then, since $C:\Lambda\times\Lambda\to
\C^\times$ is bimultiplicative, the identity \fl(comm) will continue
to hold for the bimultiplicative extension of $\e$ to the whole
$\Lambda$. Note that for $p=1$ or 2 we have 
$C(\alpha, \beta) = (-1)^{(\alpha|\alpha)(\beta|\beta)}
(-1)^{(\alpha|\beta)}$.

\begin{Thm}\label{thm:lattice}\sl
Under the assumptions {\rm (i)} and {\rm (ii)} above, 
the vertex operators $X_\alpha, X_\beta \in \F\{V\}$ are local of
order 
$$
N(\alpha, \beta) = \max \set{-m_s}{m_s<0, \ 0\le s\le p-1}\cup\{0\}.
$$
They generate the lattice vertex superalgebra 
$\goth V = \goth V_\Lambda \subset  \F\{V\}$, which does not depend on
the $\Gamma$-grading of  $\goth h$. The products of the
generators are given by 
\begin{equation}\label{fl:voprod2}
X_\alpha \ensquare{-(\alpha|\beta)-n-1} X_\beta = 
\varkappa(\alpha, \beta)\, 
%\sum_{\boldsymbol r\in \cal P(-(\alpha|\beta)-n-1)} \ 
%\prod_{j\ge1} \(\frac{\alpha(-j)}{j!}\)^{r_j} X_{\alpha+\beta},
\frac 1{n!} \big(D-\beta(-1)\big)^{(n)} X_{\alpha+\beta},
\end{equation}
for $n\ge 0$, and $X_\alpha \ensquare n X_\beta = 0$ 
if $n\ge -(\alpha|\beta)$. Here
\begin{equation}\label{fl:kappa}
\varkappa(\alpha, \beta) = \e(\alpha,\beta)\,
p^{-(\alpha|\beta)} \prod_{s=1}^{p-1}(1-\omega^s)^{m_s}.
\end{equation}
\end{Thm}
In particular,
\begin{equation*}
X_\alpha \ensquare{-(\alpha|\beta)-1} X_\beta =
\varkappa(\alpha,\beta) X_{\alpha+\beta}, \qquad
X_\alpha \ensquare{-(\alpha|\alpha)-2}X_{-\alpha} = 
\varkappa(\alpha,\alpha)\inv\,\~\alpha.
\end{equation*}

In the case when $p=1$ this is just the usual 
construction of lattice vertex algebras, see e.g. \cite{flm,kac2}.

The lattice vertex algebra 
$\goth V = \bigoplus_{\alpha\in \lambda} \goth V_\alpha$ is graded by
the lattice $\Lambda$. The subspace $\goth V_0$ is the Heisenberg
vertex subalgebra of $\goth V$, and the rest of $\goth V_\alpha$'s are
irreducible modules over $\goth V_0$. 

The even and odd parts of $\goth V$ are 
$$
\goth V\even = 
\bigoplus_{\substack{\alpha\in \Lambda:\\ (\alpha|\alpha)\in 2\Z}} \goth V_\alpha, 
\qquad
\goth V\odd = 
\bigoplus_{\substack{\alpha\in \Lambda:\\ (\alpha|\alpha)\in 2\Z+1}}
\goth V_\alpha.
$$
We also note that $\goth V$ is a simple vertex algebra.

Let us calculate the commutator map of the cocycle 
$\varkappa:\Lambda\times \Lambda\to \C^\times$ given by \fl(kappa): 
\begin{align*}
\varkappa(\alpha, \beta)\varkappa(\beta, \alpha)\inv  &= C(\alpha, \beta) \prod_{s=1}^{p-1} (1-\omega^s)^{m_s-m_{p-s}}\\
& =(-1)^{(\alpha|\alpha)(\beta|\beta)+(\alpha|\beta)} 
\prod_{s=1}^{p-1} (-\omega^{-s})^{m_s} 
\prod_{s=1}^{p-1}\(\frac{1-\omega^s}{1-\omega^{-s}}\)^{m_s}\\
& = (-1)^{(\alpha|\alpha)(\beta|\beta)} (-1)^{(\alpha|\beta)}.
\end{align*}
It follows that for different automorphisms $\sigma$
all resulting cocycles $\varkappa$ are
cohomological and define certain class 
$\varkappa \in H^2(\Lambda, \C^\times)$.
Therefore, the vertex algebra $\goth V_\Lambda\subset \F\{V\}$ 
generated by $X_\alpha$'s is indeed independent on $\sigma$.  

\begin{Rem}
If $(\alpha|\beta)$ is not an integer, then the vertex operators
$X_\alpha$ and $X_\beta$ are not local. However, they are local in a
generalized sense, and they generate a generalized vertex algebra
\cite{dl}, mentioned at the end of \sec{modules}. One can define
products $X_\alpha \ensquare n X_\beta$ for 
$n\equiv -(\alpha|\beta)\mod \Z$.
\end{Rem}

\subsection{Proof of \thm{lattice}}\label{sec:proof}
The assertion about  locality follows from
\prop{VO}\ref{VO:XX} and the formula 
$$
(w-z)=\prod_{s=0}^{p-1}\big(w^{\frac 1p}- \omega^s z^{\frac 1p}\big).
$$

Denote by $\cal P(m) = 
\bigset{\boldsymbol r=(r_1,r_2,\ldots)}{r_i\ge 0, \ \sum_{i\ge1} i r_i=m}$
the set of partitions of $m\in\Z$. Some standard combinatorial
argument shows that  \fl(voprod2) can be rewritten as

\begin{equation}\label{fl:voprod}
X_\alpha \ensquare n X_\beta = 
\varkappa(\alpha, \beta)\!
\sum_{\boldsymbol r\in \cal P(-(\alpha|\beta)-n-1)} \ 
\prod_{j\ge1} \(\frac{\alpha(-j)}{j!}\)^{r_j} X_{\alpha+\beta},
\end{equation}
where $n < - (\alpha|\beta)$. 

Consider the
operator $\delta = \imath_{w,z}-\imath_{z,w}:\C(w,z)\to
\C[[w^{\pm1},z^{\pm1}]]$, see  \fl(imath). It is easy to see
that $\delta(g)=0$ if and only if $g\in\C[w^{\pm1},z^{\pm1}]$. In this
case $g$ commutes with $\delta$: $\delta(gf)=g\delta(f)$.
We will also  make use of the
following formula. For any formal power series $f(w,z)$ in the
variables $w,z$ one has
\begin{equation*}%\label{fl:delta}
\op{Res}_w \Bigl( f(w,z)\, \delta (w-z)^{-k-1} \Bigr) = 
\frac 1{k!}\, \frac{\partial^k}{\partial w^k}\bigg|_{w=z} f(w,z)
\end{equation*}
whenever both sides make sense.

Let $n=-(\alpha|\beta)-k-1$. We calculate the product 
$X_\alpha\ensquare n X_\beta$ using \fl(lprod) and
\prop{VO}\ref{VO:XX}:  
\vskip-20pt
\begin{multline}\label{fl:XanXb}
\big(X_\alpha \ensquare n X_\beta\big)(z)= \op{Res}_w \Bigg( \Big(
X_\alpha(w) X_\beta(z)\, \imath_{w,z}(w-z)^n \\[-2pt] 
\shoveright{ -(-1)^{(\alpha|\alpha)(\beta|\beta)}
 X_\beta(z)X_\alpha(w)\, \imath_{z,w}(w-z)^n
\Big)}\\
\shoveright{\times F(N+m+k) 
\prod_{s=1}^{p-1}\big(w^{\frac 1p}- \omega^s z^{\frac 1p}\big)^{N-n}\Bigg)}\\[-5pt]
\shoveleft{=\e(\alpha, \beta) \, \op{Res}_w \Bigg(
X_{\alpha, \beta}(w,z)\, F(N+m+k)} \\[-3pt]
\shoveright{\times\prod_{s=1}^{p-1}\big(w^{\frac 1p}- \omega^s z^{\frac 1p}\big)^{N-n}
\ \delta\prod_{s=0}^{p-1}\big(w^{\frac 1p}- \omega^s z^{\frac 1p}\big)^{m_s+n}\Bigg)}\\[5pt]
\shoveleft{= \e(\alpha, \beta) \op{Res}_{w^{\frac 1p}} \Bigg(
 X_{\alpha, \beta}(w,z)\, w^{\frac{p-1}p} F(N+m+k)} \\[-6pt]
\shoveright{\times\prod_{s=1}^{p-1}\big(w^{\frac 1p}- \omega^s z^{\frac 1p}\big)^{m_s+N}\,
\delta \big(w^{\frac 1p}- z^{\frac 1p}\big)^{-k-1}\Bigg)}\\
\shoveleft{ =  \e(\alpha, \beta) \frac 1{k!}\, 
\frac{\partial^k}{\big(\partial w^{\frac 1p}\big)^k}\Bigg|_{w^{\frac 1p}=z^{\frac 1p}}}\\
\(X_{\alpha, \beta}(w,z)\, w^{\frac{p-1}p} F(N+m+k) 
\prod_{s=1}^{p-1}\big(w^{\frac 1p}- \omega^s z^{\frac
1p}\big)^{m_s+N}\).
\end{multline}
We use here that $\op{Res}_w = \op{Res}_{w^{\frac 1p}}
w^{\frac{p-1}p}$. Note that $N+m_s \ge 0$ for all $0\le s\le p-1$.

Set
$$
B=w^{\frac{p-1}p} F(N+m+k) 
\prod_{s=1}^{p-1}\big(w^{\frac 1p}- \omega^s z^{\frac 1p}\big)^{m_s+N}.
$$
Recall that there is an operator $D:V\to V$ given by \fl(Dvir), such
that the weights of all vertex operators $X_\alpha$ are 0, see
\sec{D}.
Using \fl(F) and the formula $\prod_{s=1}^{p-1}\big(1-\omega^s\big)=p$, we calculate
$$
B\,\big|_{w^{\frac 1p}=z^{\frac 1p}} = 
p^{-m-k} \prod_{s=1}^{p-1} \big(1-\omega^s\big)^{m_s}\, 
z^{-(\alpha'|\beta')}z^{k\frac{1-p}p},
$$
hence $\op{wt} B\,\big|_{w^{\frac 1p}=z^{\frac 1p}} =
-(\alpha'|\beta')-k\frac{p-1}p$. Since $B$ is a Laurent polynomial in 
$w ^{\frac 1p}$ and $z^{\frac 1p}$, we get that 
$$
\op{wt}
\frac{\partial^i}{\big(\partial w^{\frac 1p}\big)^i}
\Bigg|_{w^{\frac1p}=z^{\frac 1p}} B =
-(\alpha'|\beta')-k\,\frac{p-1}p - \frac ip.
$$

Let us now calculate the derivative of $X_{\alpha,\beta}$, as in
the proof of \prop{VO}\ref{VO:D}:
\begin{align*}
\frac{\partial}{\partial w^{\frac 1p}} X_{\alpha, \beta}(w,z) &= 
p\, w^{\frac{p-1}p}\frac{\partial}{\partial w} F_{\alpha, \beta}(w,z)\\
& = p\, w^{\frac{p-1}p}\Big(\:\~\alpha X_{\alpha, \beta}\: 
- \frac 12 (\alpha'|\alpha')\,w\inv X_{\alpha, \beta} \Big) \\
&= p\, w^{\frac{p-1}p}\Big(\~\alpha \ensquare{-1}X_{\alpha, \beta} 
+\frac12 (\alpha'|\beta')\, w\inv X_{\alpha, \beta}  \Big).
\end{align*}
Iterating this formula and using the fact that $X_{\alpha, \beta}(z,z) = 
X_{\alpha+\beta}(z)\, z^{(\alpha'|\beta')}$, we get
\begin{align}\label{fl:dX}
&\frac 1{i!}\, 
\frac{\partial^i}{\big(\partial w^{\frac 1p}\big)^i}
\Bigg|_{w^{\frac1p}=z^{\frac 1p}} X_{\alpha, \beta}\\
&\hbox to0.8\columnwidth{\hfill $= p^i\,z^{(\alpha'|\beta')} z^{i\,\frac{p-1}p}
\sum_{(r_1,r_2,\ldots)\in \cal P(i)} \ 
\prod_{j\ge1} \(\frac{\alpha(-j)}{j!}\)^{r_j} X_{\alpha+\beta}$}\notag \\ 
&\hbox to0.8\columnwidth{\hfill $+\ \ \text{a sum of terms of the
weight less than} \ i\,\frac{p-1}p + (\alpha'|\beta').$}\notag
\end{align}
Now expand \fl(XanXb) using the Leibniz rule and \fl(dX), and note that
the only term of weight 0 in this expansion is the  biggest
weight term in 
$$
\e(\alpha, \beta) \(\frac 1{i!}\, 
\frac{\partial^i}{\big(\partial w^{\frac 1p}\big)^i}
\Bigg|_{w^{\frac1p}=z^{\frac 1p}} X_{\alpha, \beta}\) B,
$$ 
which is precisely the right-hand side of the formula \fl(voprod).
This finishes the proof of \thm{lattice}.\qed

\begin{Rem}
One could have used  the theorem of Li and Xu \cite{lixu} on
characterization of lattice vertex algebras to show that the vertex
algebra generated by the vertex operators $X_\alpha\in \F\{V\}$ is
isomorphic to the lattice vertex algebra. They consider only the 
case of an even lattice, but one can easily generalize their result to the
case of a non-trivial odd part. It is easy then to deduce the product
formula \fl(voprod) but this proof would not give the explicit formula
\fl(kappa) for the cocycle $\varkappa$.

Dong and Lepowsky \cite{dl2} have constructed all twisted modules of
$\goth V_\Lambda$ for the case when $\Lambda$ is even, see
\sec{modules} for the definition. They used, however, completely different techniques, in
particular, they didn't get the cocycle $\varkappa$ explicitly.
Some special cases were also considered in \cite{dongtwisted,donag}.
We will study twisted modules of $\goth V_\Lambda$ in \sec{auto}.
\end{Rem}

\begin{Cor}\label{cor:submodule}\sl
Let $U\subset V$ be a $U(H)\,\~\otimes\, R$-submodule of $V$. Then $U$ 
is a submodule over the vertex algebra $\goth V_\Lambda$.
\end{Cor}

\subsection{Do we get all generalized modules of $\goth V_\Lambda$?}
\label{sec:lepwil}
In this section we show that in fact the construction of
\thm{lattice} exhausts all reasonable 
modules of $\goth V_\Lambda$. A similar argument was also used in 
\cite{dong,dongtwisted,lepwil81,lixu}.

Let $V$ be a generalized module over
the lattice vertex algebra $\goth V_\Lambda$ which is twisted as a
module over the Heisenberg algebra $\goth V_0$. This means that for any
homogeneous $h\in\goth h^{[\lambda]}$ the corresponding vertex operator 
$\~h \in \F\{V\}^{[\lambda]}$ is also homogeneous. 
Then the  $V$ is a module over the twisted Heisenberg algebra $H$, 
such that $\c$ acts as the identity, see \sec{heisen}. 

Assume that $V=\bigoplus_{\xi\in(\goth h^{[0]})^*}V_\xi$ is $\goth
h^{[0]}$-diagonalizable, though as in the remark preceding
\prop{VO}, we note that this assumption can be relaxed. 
We will show that $\goth V_\Lambda$-module $V$ can be obtained
by the construction of \sec{lattice}.

For $\alpha \in \Lambda$ define 
$$
e(\alpha) = E_-(\alpha, z)\inv X_\alpha(z) E_+(\alpha, z)\inv z^{-\alpha(0)}
z^{(\alpha'|\alpha')/2}.
$$
Using the same calculations as in the proof of \ref{VO:ad} and 
\ref{VO:D} of \prop{VO} we see that $\frac d{dz} e(\alpha) = 0$ and 
$\ad{h(n)}{e(\alpha)} = \delta_{n,0}\, (h|\alpha)\, e(\alpha)$, so
that  the relations \fl(RH) hold. Now we want to show that
$e(\alpha):V\to V$ generate an action of the central
extension $\widehat\Lambda$ of the lattice $\Lambda$ corresponding to
a cocycle $\e:\Lambda\times\Lambda\to \C^\times$ satisfying
\fl(comm). It will follow that $V$ is a module over the skew tensor 
product $U(H)\,\~\otimes\, R$ introduced in \sec{fock}.

Modifying slightly the proof of \prop{VO}\ref{VO:XX} we get that 
$$
X_\alpha(w)X_\beta(z) = e(\alpha) e(\beta)\, X'_{\alpha, \beta}(w,z)\,
\prod_{s=0}^{p-1}\big(w^{\frac 1p} - \omega^s z^{\frac
1p}\big)^{(\sigma^{-s}\alpha|\beta)}, 
$$
where 
\begin{multline*}
X'_{\alpha, \beta}(w,z) = 
E_-(\alpha, w) E_-(\beta, z) E_+(\alpha, w)E_+(\beta, z)\\ 
\times w^{\alpha(0)}z^{\beta(0)} w^{-(\alpha'|\alpha')/2} z^{-(\beta'|\beta')/2}
\end{multline*}
is symmetric, $X'_{\alpha, \beta}(w,z) = X'_{\beta, \alpha}(z,w)$, and also 
$e(\alpha+\beta)X'_{\alpha,\beta}(z,z) = 
X_{\alpha+\beta}(z) \, z^{(\alpha'|\beta')}$. 

Take $N$ sufficiently large. Then, since $X_\alpha$ and
$X_\beta$ must be local, we get 
\begin{align*}
0 &= \Big(X_\alpha(w) X_\beta(z)\,  - 
(-1)^{(\alpha|\alpha)(\beta|\beta)}\,
X_\beta(z) X_\alpha(w)\Big)\,  (w-z)^N\\
&= \Big(e(\alpha)e(\beta)-C(\alpha, \beta)\,e(\beta)e(\alpha) \Big)
X'_{\alpha, \beta}(w,z) 
\prod_{s=0}^{p-1}\big(w^{\frac 1p} - \omega^s z^{\frac 1p}\big)^{n_s}
\end{align*}
for some $n_s\in \Z_+$ and $C(\alpha, \beta)$ as in \fl(comm). 
But this can only happen if 
$$
e(\alpha)e(\beta) = C(\alpha, \beta)\, e(\beta)e(\alpha),
$$ 
because otherwise we get $X'_{\alpha, \beta}(w,z)\,(w-z)^n=0$ for 
$n=\max_s n_s$,
and this is impossible since $X'_{\alpha,\beta}$ is regular at $w=z$. 
 
Now we can apply the calculations \fl(XanXb) to the product 
$X_\alpha \ensquare n X_\beta$ for $n=-(\alpha|\beta)-1$. 
Similar to \fl(XanXb) we get that 
$$
e(\alpha+\beta)\, X_\alpha \ensquare n X_\beta =
e(\alpha)e(\beta)\, p^{-(\alpha|\beta)}\prod_{s=1}^{p-1}\big(1-\omega^s\big)^{(\sigma^{-s}\alpha|\beta)} \,X_{\alpha+\beta}.
$$ 
But we know that 
$X_\alpha \ensquare n X_\beta = \varkappa(\alpha, \beta) X_{\alpha+\beta}$
in $\goth V_\Lambda$, therefore 
$$
e(\alpha)e(\beta) X_{\alpha+\beta} = 
\e(\alpha, \beta) e(\alpha+\beta) X_{\alpha+\beta},
$$
hence $e(\alpha)e(\beta)=\e(\alpha, \beta)e(\alpha+\beta)$.

So we have proved the following

\begin{Thm}\label{thm:exhaust}\sl
Let $V$ be a generalized module over the lattice vertex superalgebra
$\goth V_\Lambda$
such that $\~h \in \F\{V\}^{[\lambda]}$ for $h\in \goth
h^{[\lambda]}$. Assume that $V$ is $\goth
h^{[0]}$-diagonalizable. Then there is a unique action of the extended
lattice $\^\Lambda$ on $V$ such
that for every $\alpha\in\Lambda$ the highest weights vectors
$X_\alpha\in\goth V_\alpha$  act on $V$ by the vertex
operators \fl(VO).
\end{Thm}

It follows from \thm{exhaust} that a $\goth V_\Lambda$-submodule
$U\subset V$ is a $U(H)\,\~\otimes\, R$-submodule of $V$. Combining this with
\cor{submodule}, we get the following statement.
\begin{Cor}\label{cor:irr}\sl
Assume $V$ is an irreducible module over $U(H)\,\~\otimes\, R$, satisfying the
assumptions of \thm{lattice}. Then $V$ is an irreducible module over
the vertex algebra $\goth V_\Lambda$. Conversely, let $V$ be a module
over $\goth V_\Lambda$, satisfying the assumptions of
\thm{exhaust}. Then $V$ is an irreducible $U(H)\,\~\otimes\, R$-module.
\end{Cor}

\subsection{Twisted modules over lattice vertex superalgebras}
\label{sec:auto}
In this section we study twisted modules over the lattice vertex
algebra $\goth V=\goth V_\lambda$, see \sec{modules} for the definition.
As in \sec{tvo}, we fix a $p$-th primitive root of unity 
$\omega = \exp\big(\frac{2\pi i}p\big)$.

Assume that the automorphism $\sigma:\goth h\to \goth h$ preserves the
lattice $\Lambda$. Then $\sigma$ induces an automorphism 
$\sigma:\goth V_0\to \goth V_0$ of the Heisenberg vertex algebra, 
see \sec{heisen}. Let $\^\sigma:\goth V\to \goth V$ be an extension
of this automorphism to the whole
vertex algebra $\goth V$. It is easy to see that 
$\^\sigma X_\alpha = \f(\alpha) X_{\sigma\alpha}$ for some 1-cocycle
$\f:\Lambda\to \C^\times$ such that 
$$
d\f(\alpha, \beta)= \frac{\f(\alpha+\beta)}{\f(\alpha)\f(\beta)}=
\frac{\varkappa(\sigma\alpha,\sigma\beta)}{\varkappa(\alpha,\beta)}.
$$
Since $\sigma$ preserves the norm $(\,\cdot\,|\,\cdot\,)$ on 
$\goth h$ and 
$$
\varkappa(\alpha,\beta)\varkappa(\beta,\alpha)\inv =
(-1)^{(\alpha|\alpha)(\beta|\beta)}(-1)^{(\alpha|\beta)},
$$ 
we see that  the cocycle 
$\frac{\varkappa(\sigma\alpha,\sigma\beta)}{\varkappa(\alpha,\beta)}$
is indeed symmetric and therefore equal to $d\f_0$ for 
$$
\f_0(\alpha) = 
\(\frac{\varkappa(\sigma\alpha,\sigma\alpha)}{\varkappa(\alpha,\alpha)}\)
^{\frac12}.
$$

The formula \fl(kappa) implies that 
$$
\frac{\varkappa(\sigma\alpha,\sigma\beta)}{\varkappa(\alpha,\beta)} = 
\frac{\e(\sigma\alpha,\sigma\beta)}{\e(\alpha,\beta)}, 
$$
therefore the map $e(\alpha) \mapsto \f(\alpha)\,e(\sigma\alpha)$
defines an automorphism of the algebra $R=\C[\^\Lambda]/\<\Phi\equiv\C^\times\>$. 
And vice versa, any
such automorphism of $R$ defines an automorphism of $\goth V$. 

\begin{Rem}
In  general, when the order of $\sigma$ is $p$, 
the order of $\^\sigma:\goth V\to \goth V$ may not be
equal to $p$, in fact it may not be finite at all. However, the
extension corresponding to $\f_0$ is of order $p$. Indeed, 
\begin{align*}
\^\sigma^p X_\alpha 
&= \f_0(\alpha)\f_0(\sigma\alpha) \cdots 
\f_0(\sigma^{p-1}\alpha) X_\alpha\\
&=\(\frac{
\varkappa(\sigma\alpha,\sigma\alpha)
\varkappa(\sigma^2\alpha,\sigma^2\alpha)\cdots
\varkappa(\alpha,\alpha)
}{
\varkappa(\alpha,\alpha)
\varkappa(\sigma\alpha,\sigma\alpha)\cdots
\varkappa(\sigma^{p-1}\alpha,\sigma^{p-1}\alpha)
}\)^{\frac12} X_\alpha = X_\alpha.
\end{align*}
\end{Rem}

Now we show how the vertex algebra $\goth V$ is decomposed into a sum of
the root spaces of $\^\sigma$. Consider the action of $\^\sigma$ on a
linear span \break $\op{Span}\{\,X_\alpha, X_{\sigma\alpha}, \ldots, 
X_{\sigma^{p-1}\alpha}\,\}$ of the $\^\sigma$-orbit of $X_\alpha$. 
It is easy to see that the 
eigenvalues of $\^\sigma$ on this space are all the different roots 
\begin{equation}\label{fl:root}
\mu(\alpha) = \sqrt[\uproot{3}p]{\f(\alpha)\f(\sigma\alpha)\cdots
\f(\sigma^{p-1}\alpha)}.
\end{equation}
Fix such a root $\mu=\mu(\alpha)$. 
The eigenvector corresponding to the eigenvalue $\mu\,\omega^j$
for some $0\le j\le p-1$ is
\begin{equation}\label{fl:eigen}
Y_j =\sum_{s=0}^{p-1}\omega^{-js}k_s X_{\sigma^s\alpha},
\end{equation}
where
$k_s = \mu^{-s}\, \f(\alpha)\f(\sigma\alpha)\cdots
\f(\sigma^{s-1}\alpha)$ for $s>0$ and $k_0 = 1$.

So we deduce that $\goth V = \bigoplus_{\mu\in\C^\times} \goth V_\mu$ 
is decomposed into a direct sum of the root spaces of $\^\sigma$ such
that $\^\sigma\big|\raisebox{-3pt}{$\goth V_\mu$}=\mu$.

Let $V = \bigoplus_{\xi\in (\goth h^{[0]})^*} V_\xi$ be a generalized
$\goth h^{[0]}$-diagonalizable $\goth V$-module. Denote by  
$\Xi = \bigset{\xi\in(\goth h^{[0]})^*}{V_\xi\neq 0}$ the set of weights
of $V$.  By \thm{exhaust} there are operators $e(\alpha):V\to V$,\ 
$\alpha\in\Lambda$, that define an action of $R$ on $V$, satisfying the
commutation relations \fl(RH).

\begin{Thm}\label{thm:twisted}\sl
The $\goth V$-module $V$ is $\^\sigma$-twisted 
if and only if for every
$\alpha\in \Lambda$ there is a number 
$\mu(\alpha)\in\C^\times$, given by \fl(root),  such that
\begin{itemize}
\item[\hbox to15pt{\hfill\rm(i)}]
$\displaystyle{e(\sigma^s\alpha) = k_s\inv\, e(\alpha)}$
for every $0\le s\le p-1$ and
\vskip3pt
\item[\hbox to15pt{\hfill\rm(ii)}]
$\displaystyle{\xi\big(\alpha(0)\big) \equiv 
\frac{(\alpha'|\alpha')}2 - \lambda \mod \Z}$ for every $\xi\in\Xi$,
\end{itemize}
where $p$ is the length of $\sigma$-orbit of $\alpha$, \ 
$k_s$ are given by \fl(eigen) and 
$\lambda =\frac1{2\pi i}\ln \mu(\alpha)$.
\end{Thm}

The condition (i)  means that 
$e(\sigma\alpha) = \mu(\alpha)\,\f(\alpha)\inv\, e(\alpha)$ for
every $\alpha\in \Lambda$, where the root $\mu(\alpha)$ is the same for all
$\alpha$'s in a same $\sigma$-orbit.
We also note that it is  enough to impose (i) and (ii) only
for $\alpha$ running over an integer basis of $\Lambda$.

\begin{proof}
Fix a root $\mu=\mu(\alpha)$ satisfying \fl(root), 
and let $\lambda =\frac1{2\pi i}\ln \mu$.
Let $\pi:\goth V\to \F\{V\}$ be the representation map. Then $V$ is a twisted
$\goth V$-module if and only if $\pi(\goth V_\mu) \subset 
\F\{V\}^{[\lambda]}$. It is enough to require that 
$Y_j\in \F\{V\}^{[(j/p)+\lambda]}$ 
where $Y_j$ is (the image under $\pi$ of) 
the eigenvector of $\^\sigma$ given by \fl(eigen). 
Moreover, it is enough to require this only for a finite set of
generators of $\goth V$, for example for all the $Y_j$'s corresponding
to the $\sigma$-orbits of an integer basis of $\Lambda$.

Set $\overset{\circ}{X}_\alpha(z)=
E_-(\alpha,z)E_+(\alpha,z)\subset \F\{V\}$  so that (see \sec{tvo})
$$
X_\alpha = e(\alpha)\,\overset{\circ}{X}_\alpha\,
z^{\alpha(0)}z^{-(\alpha'|\alpha')/2}.
$$ 
Set also $\overset{\circ}{Y}_j=
\sum_{s=0}^{p-1}\omega^{-js}\,k_s\,
e(\sigma^s\alpha)X_{\sigma^s(\alpha)}$ and then 
$$
Y_j = \overset{\circ}{Y}_j\,
z^{\alpha(0)}z^{-(\alpha'|\alpha')/2}.
$$  
It follows that the field $Y_j \in \F\{V\}$ is homogeneous if and only if
$\overset{\circ}{Y}_j$ is homogeneous and the values
$\xi\big(\alpha(0)\big)$  
are the same modulo $\Z$ for all $\xi\in\Xi$.

Assume that $\overset{\circ}{Y}_j$ is homogeneous.
Since $\overset{\circ}{X}_\alpha(z)\in \bigoplus_{q=0}^{p-1} \F\{V\}^{[q/p]}$,
we have
$\overset{\circ}{Y}_j\in \F\{V\}^{[(q+j)/p]}$
for some $0\le q \le p-1$. 
Let $\tau$ be the automorphism
of $\bigoplus_{q=0}^{p-1} \F\{V\}^{[q/p]}$ such that 
$\tau\big|\raisebox{-3pt}{$\F\{V\}^{[q/p]}$} = \omega^q$. It is easy
to see that $\overset{\circ}{X}_{\sigma\alpha} = 
\tau\, \overset{\circ}{X}_\alpha$. 
Take some $m\equiv -\frac rp \mod \Z$. Denote by $x_s$ the coefficient of
$z^m$ in $\overset{\circ}{X}_{\sigma^s\alpha}$. It follows that 
$x_s = \omega^{rs} x_0$. If $r\not\equiv q+j \mod p$, then the
coefficient of $z^m$ in $\overset{\circ}{Y}_j$ is equal to 0. This
gives us the following system of linear equations:
$$
\sum_{s=0}^{p-1}\omega^{(r-j)s}\, k_s\, e(\sigma^s\alpha) =0 \quad
\text{for}\quad 0\le r\le p-1, \ \ r\not\equiv q+j \mod p.
$$
Since $k_0$=1, the solution of this system is $e(\sigma^s\alpha) = 
\omega^{qs}\,k_s\inv \,e(\alpha)$. Taking
$\mu\,\omega^q$ instead of $\mu$ we get exactly (i). 

The condition (ii) follows from the fact that in order to have 
$\deg Y_j = \frac jp+\lambda$ we must  have  
$$
\xi\big(\alpha(0)\big) \equiv 
\frac{(\alpha'|\alpha')}2 -\lambda-\frac qp  \equiv 
\frac{(\alpha'|\alpha')}2-\frac1{2\pi i}\ln \big(\mu\,\omega^q\big)
\mod \Z
$$ 
for all $\xi\in\Xi$.
\end{proof}

It follows that sometimes the lattice vertex algebra 
$\goth V=\goth V_\Lambda$ 
has no twisted representations that agree with an
automorphism $\sigma:\Lambda\to\Lambda$.
Recall that $C:\Lambda\times\Lambda \to \C^\times$ is the
commutator map given by \fl(comm).

\begin{Cor}\label{cor:zero}\sl
If for some  $\alpha\in \Lambda$ and $0\le j\le p-1$ we have 
$C(\alpha, \sigma^j\alpha)\neq 1$, then there are no non-trivial
twisted modules of $\goth V_\Lambda$ corresponding to either of the
extensions $\^\sigma:\goth V_\Lambda\to\goth V_\Lambda$ of $\sigma$. 
\end{Cor}
\begin{proof}
Set $\beta =\sigma^j\alpha$.
If $C(\alpha,\beta)\neq1$, then  $e(\alpha)e(\beta)\neq e(\beta)e(\alpha)$.
Yet by \thm{twisted}, we have
$e(\alpha)e(\beta)\inv \in \C^\times$. Therefore, $R$ acts by 0 on every
twisted $\goth V_\Lambda$-module, 
hence by \thm{exhaust}, so does $\goth V_\Lambda$. 
\end{proof}

\subsection{Semisimplicity of twisted representations}
\label{sec:ss}
Consider the category $\cal O_{\^\sigma}$ of 
$\goth h^{[0]}$-diagonalizable $\^\sigma$-twisted 
 modules $V$ over $\goth V_\Lambda$ such that 
$V$ satisfies the conditions of \lem{heisen} as a module over the
Heisenberg algebra $H$. In particular, $V$ must be a module over the
vertex operator algebra $\goth V_\Lambda$. Note that if 
$V=\bigoplus_n V_n$ is a $\^\sigma$-twisted module
over the vertex operator algebra $\goth V_\Lambda$ such
that $V_n=0$ for $n\ll0$, then $V\in \cal O_{\^\sigma}$. 
The latter modules
appear in the representation theory of vertex algebras quite often,
in particular in connection with the Zhu theory \cite{dlm,zhu}.   

Note that by \cor{zero}, the category  $\cal O_{\^\sigma}$ can
sometimes be trivial. 

\begin{Thm}\label{thm:ss}\sl
The category $\cal O_{\^\sigma}$  is semisimple
with finitely many isomorphism classes of simple objects.
\end{Thm}

In order to prove this theorem we have to study 
the quotient algebra $A$ of $R=\C[\^\Lambda]/\<\Phi\equiv\C^\times\>$ 
modulo the ideal generated by
relations (i) of \thm{twisted}.

More precisely, 
let $\Pi\subset\Lambda$ be a finite set of vectors, closed under
$\sigma$ and spanning $\Lambda$ over $\Z$. For $\alpha\in\Pi$
choose the roots $\mu(\alpha)$, given by \fl(root), such that 
$\mu(\alpha) = \mu(\beta)$ if $\alpha$ and $\beta$ lie in the same
$\sigma$-orbit of $\Pi$.
For $x\in R$ denote by $\ol x$ its image in $A$. 
Then the algebra $A$ is generated by the set  
$\bigset{\ol{e(\alpha)}}{\alpha\in\Pi}$ subject to relations 
$\ol{e(\sigma \alpha)} = \mu(\alpha)\,\phi(\alpha)\inv\,\ol{e(\alpha)}$ for
all $\alpha\in \Pi$. Hence $A$ depends on the
cocycle $\f:\Lambda\to \C^\times$, which determines the extension
$\^\sigma$, and also on the choice of roots
$\mu(\alpha)$ for every $\sigma$-orbit of $\Pi$.

Note that the relations (i) of \thm{twisted} belong in fact to
the group $\^\Lambda$. Let $G$ be the quotient group of  $\^\Lambda$
modulo the normal subgroup generated by these relations. Then 
$A=\C[G]/\<\Phi\equiv\C^\times\>$.

Let $\nu:\Lambda\to \big(\goth h^{[0]}\big)^*$ be the map given by 
$\nu(\alpha)h = (\alpha|h)$. Then $\nu(\Lambda)$ is a sublattice of
the dual lattice to $\Lambda^{[0]}$. 
The algebra $R$ is graded by 
$\nu(\Lambda)$ by setting $\deg e(\alpha)=\nu(\alpha)$. 
We observe that if $\alpha,\beta \in \Lambda$ belong to the same
$\sigma$-orbit, then $\nu(\alpha) = \nu(\beta)$. It follows that the
relations (i) of \thm{twisted} are in the kernel of the
composition map $\^\Lambda \to \Lambda \xrightarrow{\nu} \nu(\Lambda)$,
which therefore induces a homomorphism $G\to \nu(\Lambda)$, and
this makes $\C[G]$ and $A$  graded by $\nu(\Lambda)$ as well. 

\begin{Prop}\label{prop:ss}\sl
The algebra $A$ is $\nu(\Lambda)$-graded semisimple.
\end{Prop}

Recall that a graded algebra is called graded semisimple if it is semisimple
as a left graded module over itself. Like in a non-graded case, an
algebra is graded semisimple if and only if it is decomposed into a
direct product of graded simple algebras, the latter by a graded
version of the classical Wedderburn theorem are isomorphic to matrix
algebras over graded division algebras. If an algebra $A$ is graded 
semisimple, then the category of graded $R$-modules is semisimple,
i.e. every module is completely reducible. The simple graded
$A$-modules are the simple homogeneous ideals of $A$ with a possible
shift of degrees. 

Similar to the non-graded case, the graded Jacobson radical 
$J_{\text{gr}}(A)$ of a graded algebra $A$ is defined to be the
intersection of all maximal graded left (or right) ideals. It is a
standard excercise to prove that  $J_{\text{gr}}(A)$ is in fact 
a double sided ideal that acts by 0 on any simple graded $A$-module. 
If $J_{\text{gr}}(A)=0$ that $A$ is called semiprimitive.

The algebra  $A$ is called graded (left) Artinian if
all strictly decreasing chains of graded left ideals are finite. In this
case $J_{\text{gr}}(A)$ is an intersecton of finitely many left
maximal ideas. If $A$ is graded Artinian and $J_{\text{gr}}(A)=0$
then $A$ is graded semisimple. 

For the non-graded case all this is a classical theory 
presented in most graduate algebra textbooks, like Jacobson's {\it
Basic Algebra} \cite{jacobson}. For the graded case the references are much
more scarce, see, however, \cite{NvO}. 

We will use the following
rather obvious fact: A homomorphic image of a graded semisimple algebra
modulo a homogeneous ideal is again a graded semisimple algebra.

\begin{proof}[Proof of \prop{ss}]
Let $\Pi =
\bigsqcup_{j=1}^n \Pi_j$ be the decomposition of $\Pi$ into a 
disjoint union of $\sigma$-orbits.  Choose some $\alpha_j\in \Pi_j$ for
every $1\le j\le n$. Let $x_j=\ol{e(\alpha_j)}\in A$.
Then the set $x_1^{\pm1},\ldots,x_n^{\pm1}$ is a set of generators of
$A$. We note that $x_i x_j = c_{ij}\,x_j x_i$ for  
$c_{ij}=C(\alpha_i,\alpha_j)\in \C^\times$, where $C$ is the
commutator map given by \fl(comm).

Let $\xi_j = \deg \alpha_j = \nu(\alpha_j)\in \nu(\Lambda)$. It is
easy to see that $\xi_j = 0$ if and only if 
$\sum_{\alpha\in \Pi_j}\alpha = 0$. We can assume that $\xi_1, \ldots,
\xi_m\neq0$ for some $m\le n$ and $\xi_{m+1}=\ldots=\xi_n=0$.
Clearly, $\xi_j$'s span $\nu(\Lambda)$ over $\Z$. 

Assume also that our choice of $\Pi$ yields the minimal possible value
of $m$. We claim that in this case $\xi_1, \ldots \xi_m$ is a
$\Z$-basis of $\nu(\Lambda)$. Indeed, otherwise there is an invertible
matrix $\sf M \in SL(m,\Z)$ such that the $m$-th column of 
$(\xi_1 \ \cdots \ \xi_m)\,\sf M$ is 0. Let $\alpha_j'$ be the
$j$-th column of $(\alpha_1 \ \cdots \ \alpha_m)\,\sf M$ for 
$1\le j\le m$, and let $\Pi'$ be the closure of 
$\left\{\alpha_1',\ldots,\alpha_m'\right\}\cup \Pi_{m+1}\cup\cdots\cup \Pi_n$
with respect to the action of $\sigma$.
Then $\Pi'$ is a generating set of $\Lambda$, closed under $\sigma$
and, since $\nu(\alpha_m')=0$ the number of $\sigma$-orbits in $\Pi'$
with non-zero degree 
is at most $m-1$. 

Suppose there is
a linear relation between elements of $\Pi$  of the form 
$$
\sum_{j=1}^n\sum_{s=0}^{p-1} r_{js}\,\sigma^s\alpha_j=0,\qquad
r_{js}\in\Z.
$$
In $R$ this relation becomes
$\prod_{j,s} e(\sigma^s\alpha_j)^{r_{js}} = \theta'$,
where $\theta' \in \C^\times$ is a product of some values of the
cocycle $\e$. In $A$ this relation becomes
\vskip-1pt
$$
\prod_{j=1}^n x_j^{r_j}=\theta,\qquad
r_j = \sum_{s=0}^{p-1}r_{js},
$$
for some other constant $\theta \in \C^\times$.

For $m+1\le j\le n$ we have $\sum_s \sigma^s\alpha_j =0$, hence in $A$
we get the relation $x_j^p=\theta_j\in \C^\times$. By the change of
variables $x_j\mapsto \theta_j^{-1/p}\,x_j$ we can make $\theta_j=1$.

A look at the formula \fl(comm) for the commutator map $C$ shows that if 
$\sum_{s=0}^p \sigma^s\alpha=0$ for some $\alpha\in\Lambda$, 
then $C(\alpha,\beta)^p = 1$ for any other $\beta \in \Lambda$.

Summing up, we get that $A$ is a
homomorphic image of the algebra 
\begin{equation*}%\label{fl:B}
B=\Bbbk\,\Bigg\langle\,x_1^{\pm1},\ldots,x_n^{\pm1}\,\Bigg|\
\begin{gathered}
x_i x_j = c_{ij}\,x_j x_i  \ \
\text{for} \ \ 1\le i,j\le n\\ 
x_j^p = 1 \ \ \text{for}\ \ m+1\le j\le n
\end{gathered}
\,\Bigg\rangle,
\end{equation*}
where $c_{ij}=c_{ji}\inv\in\C^\times$ for $1\le i,j\le n$ and 
$c_{ij}^p=1$ if either $i>m$ or $j>m$. The $\nu(\Lambda)$-grading on
$B$ is defined by $\deg x_j = 0$ for 
$m+1\le j\le n$ and $\deg x_j =\xi_j$ for $1\le j \le m$ where 
$\xi_1,\ldots,\xi_m$ is a $\Z$-basis of $\nu(\Lambda)$.
It is enough to show that $B$ is
$\nu(\Lambda)$-graded semisimple.

Let $\xi = k_1\xi_1 +\ldots +k_m \xi_m, 
\ k_i \in \Z$, be an arbitrary vector in $\nu(\Lambda)$. Denote
by $B_\xi$ the homogeneous component of $B$ of degree $\xi$. Note that 
the element $x^\xi = x_1^{k_1}\cdots x_m^{k_m} \in B_\xi$ is
invertible. For a graded subspace $I \subset B$ denote by $I_\xi$ the
homogeneous component of $I$ of degree $\xi$.

The component $B_0$ is the homomorphic image of the group
algebra of finite group
$$
G_0 = \Bigg\langle x_k, c_{ij}\,\Bigg|\,
\begin{gathered}
m+1\le i,j,k\le n,\ x_k^p = c_{ij}^p=1,\ c_{ij} = c_{ji}\inv,\\
x_ix_j = c_{ij}\,x_jx_i, \ x_k c_{ij} = c_{ij} x_k
\end{gathered}
\Bigg\rangle.
$$
By Maschke's theorem (see e.g. \cite{jacobson}), $\C[G_0]$ and hence
$B_0$ is a finite dimensional semisimple algebra. 

We claim that if $I,J \subset B$ are two graded left ideals such that
$I_0 = J_0$, then $I=J$. Indeed, if $b \in I_\xi$, then $x^{-xi}b \in
I_0 = J_0$, hence $b = x^\xi x^{-\xi} b \in J_\xi$.
It follows that $B$ is graded (left) Artinian.
 
Let $I\subset B$ be a graded left ideal. Then $I$ is graded maximal if
and only if $I_0$ is a maximal ideal of $B_0$. Therefore, since
$J(B_0) = 0$, we must have $J_{\text{gr}}(B) = 0$, and that finishes
the proof. 
\end{proof}

\begin{Rem}(D.~Passman)\quad
The argument above shows in fact that if $B=B_0*G$ is a crossed
product of an algebra $B_0$ with a group $G$, then $B$ is $G$-graded
Artinian if and only if $B_0$ is Artinian, and also $J_{\text{gr}}(B)
= 0$ if and only if $J(B)=0$.  The same is true if
$B$ is a strongly $G$-graded algebra. See e.g. \cite{karp,passman} for
definitions and further results.
\end{Rem}

\begin{proof}[Proof of \thm{ss}]
Let $V\in \cal O_{\^\sigma}$.
By (i) of \lem{heisen} and \thm{twisted} we have 
$V = M(1)\otimes \Omega$, where $\Omega$ is $\nu(\Lambda)$-graded 
module over $A$. Hence it follows from \prop{ss} that $\Omega$ is
decomposed into a direct sum of graded irreducible $A$-modules, and
\cor{irr} implies that $\goth V$ is decomposed into a direct sum of
irreducible $\goth V_\Lambda$-modules. 

A $\goth V$-module $V = M(1)\otimes \Omega\in\cal O_{\^\sigma}$ is
simple if and 
only if $\Omega$ is a simple $\nu(\Lambda)$-graded $A$-module. Such
$\Omega$ must be isomorphic to a simple homogeneous ideal of $A$ up to
a shift of weights. The weights of $\Omega$ are restricted by (ii) of
\thm{twisted}. We claim that a simple object of 
$\cal O_{\^\sigma}$ is  determined up to an isomorphism by a choice of 
roots $\mu(\alpha_j)$, given by \fl(root), for
each generating $\sigma$-orbit $\Pi_j \ni \alpha_j$ of $\Lambda$,
a choice of simple
homogeneous ideal $\Omega$ of $A$ and an equivalence class 
$\eta \in \big(\Lambda^{[0]}\big)'\big/\nu(\Lambda)$.

Indeed, assume all these choices are made. Since the extension
$\^\sigma$ is fixed, the choice of $\mu(\alpha_j)$'s determines the
$\nu(\Lambda)$-graded semisimple algebra $A$. 
The set $\Xi =\bigset{\xi \in \big(\goth h^{[0]}\big)^*}{\Omega_\xi \neq0}$ is
an equivalence class in $\big(\goth h^{[0]}\big)^*\big/\nu(\Lambda)$. 
Let $\xi \in \big(\goth h^{[0]}\big)^*$ be such that 
$\xi\big(\alpha(0)\big) 
= \frac 12 (\alpha'|\alpha')-\frac1{2\pi i}\ln \mu(\alpha)$ for
$\alpha\in\Lambda$.  
By \thm{twisted} (ii) we have that 
$\Xi \equiv \xi \mod \big(\Lambda^{[0]}\big)'$, so now we further
specify $\Xi \equiv \xi+\eta \mod \nu(\Lambda)$.  

It follows that there
are at most finitely many isomorphism classes of simple objects in 
$\cal O_{\^\sigma}$. 
\end{proof}

\subsection{Examples}\label{sec:examples}
\medskip\noindent{\it Example 0.}\quad
If $\sigma=\op{Id}$, then  $\nu(\Lambda)\cong\Lambda$, since the form
is non-degenerate, and $A=R$ is a
$\Lambda$-graded division algebra. 
The automorphism $\^{\op{id}}:\goth V_\Lambda\to\goth V_\Lambda$ is
defined by choosing an arbitrary values $\mu_j=\f(\alpha_j)\in \C^\times$
for a basis $\alpha_1,\ldots,\alpha_l$ of $\Lambda$. The simple
objects of  of category $\cal O_{\^{\op{id}}}$ are parametrized by
$\Lambda'/\Lambda$, which agrees with the result of Dong \cite{dong}. 
Let $\xi\in \goth h^*$ be the functional defined by 
$\xi(\alpha_j) = \frac 1{2\pi i} \ln \mu_j$. Then 
the simple object of category $\cal O_{\^{\op{id}}}$ corresponding to 
an equivalence class $\eta\in \Lambda'/\Lambda$ is $M(1)\otimes
\C[\Lambda]$, where $\C[\Lambda]$ is graded such that 
its set of weights is equal to $\eta+\xi+\Lambda$.

\medskip\noindent{\it Example 1.}\quad
A more interesting example is when $\sigma=-1$. In this case take
$\Pi = \left\{\pm\alpha_1,\ldots,\pm\alpha_l\right\}$.
The grading is trivial, since $\Lambda^{[0]}=0$. 
One can choose the cocycle 
$\e$ so that $\e(\alpha_j,\alpha_j)=1$ and then 
$\f(-\alpha_j)=\f(\alpha_j)\inv$, hence \fl(root) gives 
$\mu_j = \pm 1$ for all
$1\le j\le l$. Set $x_j = \ol{e(\alpha_j)}\in A$, then 
$$
A=\C[x_1, \ldots, x_l]\,\big/\<x_i x_j = c_{ij}\,x_j x_i, \  
x_j^2=\mu_j\phi(\alpha_j)\>,
$$
where $c_{ij}=(-1)^{(\alpha_i|\alpha_i)(\alpha_j|\alpha_j)}
(-1)^{(\alpha_i|\alpha_j)}$. 
So  $A$ is a semisimple algebra of
dimension $2^l$ for any choice of $\f$ and $\mu_j$.
A simple module $V\in \cal O_{\^\sigma}$ is decomposed
$V=M(1)\otimes\Omega$, where $\Omega$ isomorphic to a simple ideal of
$A$. This case was studied in \cite{dongtwisted,donag,flm}.

We remark that $-1$ is an isomorphism of any lattice $\Lambda$. When 
$\Lambda$ is the Leech lattice, a certain $\^{-1}$-twisted module over 
$\goth V_\Lambda$ is used to construct the
Moonshine vertex algebra $V^\natural$, see \cite{flm}.

\medskip\noindent{\it Example 2.}\quad
Take $\Lambda = \Z\alpha + \Z\beta$ and let $\sigma$ be the rotation
$\alpha\mapsto \beta, \ \beta\mapsto -\alpha$. Then we must have 
$(\alpha|\beta)=0, \ (\alpha|\alpha)=(\beta|\beta)$. As in the
previous example, $\Lambda^{[0]}=0$, so the grading is trivial. 
The formula \fl(comm) yields $C=1$, therefore the cocycle $\e$ is
trivial and hence 
$\f:\Lambda\to \C^\times$ is a character. Set $\f(\alpha)=\f_1, \
\f(\beta)=\f_2$. The lattice
$\Lambda$ is generated by a single orbit 
$\{\,\alpha, \beta, -\alpha, -\beta\,\}$, and for that orbit 
$\mu = \sqrt[4]1$ by \fl(root). Set $x=\ol{e(\alpha)},\,
y=\ol{e(\beta)}\in A$. Then relations (i) of \thm{twisted} give
$y=\mu\f_1\inv x, \ x\inv = \mu\f_2\inv y$, therefore $A = \C[x]/\<x^2
= \mu^2\f_1\f_2\>$. It follows that there are exactly 2 irreducible
$\^\sigma$-twisted $\goth V_\Lambda$-modules for every extension 
$\^\sigma:\goth V_\Lambda\to\goth V_\Lambda$.

\begin{acknowledgment}
I would like to thank Bojko Bakalov, Chongying Dong, Jim
Lepowsky and Don Passman for the inspiring discussions. I am very grateful to Bojko
Bakalov for communicating the unpublished manuscript \cite{bkt}.
\end{acknowledgment}

\bibliography{conformal,my,general,physics,vertex}

\end{article}

\end{document}